\documentclass[12pt]{article}
\usepackage{amsmath,amssymb,colordvi}
\usepackage{graphicx}
\usepackage{url}

\usepackage{verbatim}
\newtheorem{theorem}{Theorem}
\newtheorem{corollary}[theorem]{Corollary}
\newtheorem{lemma}[theorem]{Lemma}
\newtheorem{example}[theorem]{\it Example}
\newtheorem{proposition}[theorem]{Proposition}

\newtheorem{remark}[theorem]{\it Remark}

\newcommand{\CaixaPreta}{\vrule Depth0pt height5pt width5pt}
\newcommand{\bgproof}{\noindent {\bf Proof.} \hspace{2mm}}
\newcommand{\edproof}{\hfill \CaixaPreta \vspace{3mm}}

\def\NN{\mathbb N}
\def\CC{\mathbb C}
\def\C{\mathbb C}

\def\R{\mathbb R}
\def\RR{\mathbb R}
\def\ZZ{\mathbb Z}

\def\LL{\mathcal L}

\def\pp{\frac{1}{p}}
\def\Tpa{\mathcal{T}_p^{(\alpha)}(t^\alpha)}
\def\Tpabs{\mathcal{T}_p^{(\alpha)}(|t|^\alpha)}
\def\Sbmu{\mathcal{S}_{\beta, \mu}}

\title{Generalized Stieltjes and other integral operators on Sobolev-Lebesgue
spaces}

\author{Pedro J. Miana and Jes\'us Oliva-Maza \thanks{Authors  have been partially supported by Project MTM2016-77710-P, DGI-FEDER, of the MCYTS and Project E26-17R, D.G. Arag\'on, Spain.}
\\
\small Departamento de Matem\'aticas\\[-0.8ex]
\small Instituto Universitario de Matem\'aticas y Aplicaciones\\[-0.8ex]
\small Universidad de Zaragoza\\
\small 50009 Zaragoza, Spain\\
\small \texttt{pjmiana@unizar.es,} \texttt{joliva@unizar.es}
}

\date{}

\begin{document}
\maketitle

\begin{abstract}
For $\mu>\beta>0$, the generalized Stieltjes operators
$$
\mathcal{S}_{\beta,\mu} f(t):={t^{\mu-\beta}}\int_0^\infty
{s^{\beta-1}\over (s+t)^{\mu}}f(s)ds, \qquad t>0,
$$
 defined on
Sobolev spaces $\mathcal{T}_p^{(\alpha)}(t^\alpha)$  (where $\alpha\ge 0$ is the fractional order of derivation and these spaces are
 embedded in $L^p(\RR^+)$ for  $p\ge 1$) are studied in detail. If $0 < \beta - \pp <  \mu$, then  operators $\mathcal{S}_{\beta,\mu}$ are
  bounded (and we compute their operator norms which depend on
$p$);  commute and factorize with  generalized Ces\'{a}ro operator  on $\mathcal{T}_p^{(\alpha)}(t^\alpha)$ . We calculate and represent explicitly
  their spectrum set $\sigma (\mathcal{S}_{\beta,\mu})$. The main technique is to subordinate these operators in terms of $C_0$-groups and transfer new properties from some special functions to Stieltjes operators. We also prove some similar results for generalized Stieltjes operators
$
\mathcal{S}_{\beta,\mu}$ in the Sobolev-Lebesgue $\mathcal{T}_p^{(\alpha)}(\vert t\vert^\alpha)$ defined on the real line $\R$. We show connections with the Fourier and the Hilbert transform and  a convolution product defined by the Hilbert transform.
\end{abstract}

\date{}

%\maketitle

\section{Introduction}

%\section{Convoluted Fractional Banach Algebras on $\mathbb{R}_+$}
\setcounter{theorem}{0} \setcounter{equation}{0}

 The integral operator $\mathcal{S}$, where
$$
\mathcal{S}f(t) := \int_0^\infty \frac{f(s)}{t+s}ds, \qquad t>0,
$$
is the origin of different theories in many areas of  Mathematical Analysis and Differential Equations in real  and complex variable. For example, this operator  arises naturally from the twice iteration  of the classical Laplace transform.

Perhaps, T.J.  Stieltjes was  one of the first who studied  deeply this integral operator. He treated the distribution of a measure of density $f$ on the real positive axis and the problem of the inversion, i.e. to calculate $f$ known ${\mathcal S}f$, by
means of contour integration (\cite[p.473]{[St]}). Later D.V. Widder proved  a real inversion theorem depending only
on a knowledge of $\mathcal{S}f$ and its derivatives on the positive real axis (\cite{[Wi]}).

Given two real sequences $(a_m)_{m\ge 1} $ and $ (b_n)_{n\ge 1}$ such that $\sum_{m\ge 1} a_m^2,\,\, \sum_{n\ge 1} b_n^2<\infty$, the celebrated Hilbert's double serie inequality states that
$$
\sum_{n,m=1}^\infty{\vert a_mb_n\vert\over n+m}\le \pi\left(\sum_{m= 1}^\infty a_m^2\right)^{1\over 2} \left(\sum_{n= 1}^\infty b_n^2\right)^{1\over 2}.
$$
In \cite{Schur}, I. Schur proved
that $\pi$ is the best possible constant and also discovered the integral analogue of this inequality
which became known as the Hilbert integral inequality in the form
$$
\int_0^\infty\int_0^\infty{\vert f(t)g(s)\vert \over s+t}dtds\le \pi\left(\int_{0}^\infty f^2(t)dt\right)^{1\over 2} \left(\int_{0}^\infty g^2(t)dt\right)^{1\over 2},
$$
i.e. ${\mathcal S}:L^2(\R^+)\to L^2(\R^+)$ is a linear and bounded operator and $\Vert {\mathcal S}\Vert=\pi$.

H. Weyl (\cite{Weyl}) and T. Carleman (\cite{Ca}) continued the study of integral equations
$$
f(t)-\lambda\int_0^\infty K(s,t)f(s)ds=\varphi(t), \qquad t>0, \quad f,\varphi\in L^2(\R^+).
$$
Here the kernel $K$ is the quotient of two positive polynomials $P$ and $Q$ of degree $m$ and $m+1$ respectively; in particular the case  $K(s,t)={1\over s+t}$ is included.  This point of view was followed in \cite[Theorem 319]{Ha-Li-Po-64}, where under the same conditions about the kernel $K$, is proved that the integral operator
$$f\mapsto \int_0^\infty K(s,\cdot)f(s)ds,\qquad f\in L^p(\R^+),
$$
is bounded  on $L^p(\R^+)$ for $p>1$. In particular for $K(s,t):=\displaystyle{t^{\mu-\beta}s^{\beta -1}\over(t+s)^\mu}$, the generalized Stieltjes operator $\mathcal{S}_{\beta, \mu} $ where
\begin{eqnarray*}
\mathcal{S}_{\beta, \mu} f(t) := t^{\mu-\beta}\int_0^\infty \frac{s^{\beta -1}}{(t+s)^\mu}f(s)ds,
\end{eqnarray*}
is bounded on $L^p(\R^+)$ and
$
\Vert \mathcal{S}_{\beta, \mu}\Vert= B( \beta-1/p,\mu-\beta+1/p)
$
for $0 < \beta - \pp< \mu$ (we write ${\mathcal S_{1,1}}={\mathcal S}$).  The operator $\mathcal{S}_{\beta, \mu}$ is also bounded in $L^1(\R^+)$ under the restriction $0< \beta-1  < \mu$. Other generalization of Stieltjes operators (depends on three paremeters) are considered in \cite[Exercise (4.6.12)]{Oki}. Recently a explicit formula for the resolvent of the Carleman operator $R(z)=({\mathcal S}-zI)^{-1}$  has been calculated in \cite[Theorem 2.3]{[Yafa]}. Note also that the generalized Stieltjes transform is an iterative Laplace transform in sense that
\begin{equation}\label{faclap} \mathcal{S}_{\beta,\mu} f(t) = {1 \over \Gamma(\mu)}t^{\mu-\beta} \mathcal{L}\left(
x^{\mu-1}\mathcal{L}\left(s^{\beta-1}f(s)\right)(x)\right)(t), \qquad t>0,\,\,f\in L^p(\R^+),\end{equation}
for $0< \beta -1/p < \mu$  and $\mathcal{L}$ is the usual Laplace transform.

In this paper, our main idea  is to subordinate the generalized Stieltjes operator $\mathcal{S}_{\beta, \mu}$ in terms of a $C_0$-group of isometries $(T_{t,p})_{t\in \R}$ defined on two families of Sobolev-Lebesgue spaces $\mathcal{T}^{(\alpha)}_p(t^\alpha)$ and $\mathcal{T}^{(\alpha)}_p(\vert t\vert^\alpha)$ , which are Banach spaces continuously contained respectively in $L^p(\RR^+)$ and $L^p(\RR)$, $1\le p<\infty$,  where
\begin{equation}\label{isometry}
T_{t,p} f(s):=e^{-\frac{t}{p}}f(e^{-t}s), \qquad f \in
\mathcal{T}_p^{(\alpha)}(t^\alpha), \mathcal{T}^{(\alpha)}_p(\vert t\vert^\alpha).
\end{equation}

This strategy has been pursued by other authors. This $C_0$-group of isometries was considered in \cite{Cowen} to study the property of subnormality of the Ces\`{a}ro operator on $L^2(\R^+)$. Later in \cite{arvanti}, the Ces\`{a}ro operator on  the Hardy spaces of the half plane was considered again and finally in \cite{[LMPS]} the generalized Ces\`{a}ro operator $\mathcal{C}_\gamma$ where
$$
\mathcal{C}_\gamma f(t) := \frac{\gamma}{t^\gamma}\int_0^t (t-s)^{\gamma-1}f(s)ds,  \qquad t>0,
$$
with  $\gamma > 0$ was also considered  on the spaces  $\mathcal{T}_p^{(\alpha)}(t^\alpha)$ and $\mathcal{T}^{(\alpha)}_p(\vert t\vert^\alpha)$ for $p>1$ and $\alpha \geq 0$. Recently in  \cite{[ASV]}, authors have applied this $C_0$-group to the Stieltjes operator ${\mathcal S}$ on the  Hardy spaces of the half plane.

This subordination process can be described as an extension  of the Fourier transform in abstract Banach spaces. Let $X$ be a Banach space, $\mathcal{B}(X)$ the set of all linear and bounded operators on a
Banach space $X$, and $(T(t))_{t\in \R}\subset \mathcal{B}(X)$ a $C_0$-group of uniformly bounded operators on $X$, i.e., $T(t+s)=T(s)T(t)$, for $t, s\in \R$; $\lim_{t\to 0} T(t)x=x$ for $x\in X$ and $M:=\sup_{t\in \R} \Vert T(t)\Vert<\infty$ (\cite[Definition 3.1.19]{ABHN}). Let $\theta$ denote the  map: $\theta: L^1(\R) \longrightarrow \mathcal{B}(X)$ such that
\begin{equation}\label{subor}
\theta(g)x := \int_{-\infty}^{\infty}g(t)T(t)xdt, \qquad x\in X, \qquad g\in L^1(\R).
\end{equation}
Then, the map $\theta$ is actually well defined, is a linear and bounded operator, $\Vert \theta(g)\Vert \leq  M\|g\|_1$ for $g\in L^1(\R),$ and   $\|\theta\| \leq M$ (\cite[Lemma IV.3.17]{En-Na-00}). Also, $\theta$ is commutative in its range $\theta(L^1(\R))\subset \mathcal{B}(X)$, because equalities $\theta(f)\theta(g)= \theta (f\ast g)=\theta(g) \theta(f)$ hold for $f,g \in L^1(\R)$. Moreover, we may identify spectrum sets of bounded operators  by the spectral mapping theorem,
$$
\sigma(\theta(g))=\overline{\widehat{g}(\sigma(iA))}
$$
where $\widehat{g}$ is the Fourier transform of $g$ and $A$ is the infinitesimal generator of the $C_0$-group (see for example \cite[Theorem 3.1]{[Se]}).

The Sobolev-Lebesgue spaces
$\mathcal{T}^{(\alpha)}_p(t^\alpha)$ and $\mathcal{T}^{(\alpha)}_p(\vert t\vert^\alpha)$   were introduced and studied in detail for $p=1$ in \cite{Ga-Mi-06} and for $p>1$ in \cite{Ro-08}; see also \cite[Section 2 and 4]{[LMPS]}. In particular, $\mathcal{T}_p^{(0)}(t^0) = L^p(\R^+)$ and $\mathcal{T}_p^{(0)}(|t|^0) = L^p(\R)$, so every result obtained in this article for $\alpha \geq 0$ is also applicable to $L^p$ spaces. Moreover, the subspace
$\mathcal{T}^{(\alpha)}_p(t^\alpha)$ (and $\mathcal{T}^{(\alpha)}_p(\vert t\vert^\alpha)$) is a module for the algebra
$\mathcal{T}^{(\alpha)}_1(t^\alpha)$ (and for $\mathcal{T}^{(\alpha)}_1(\vert t\vert^\alpha)$) for the convolution product
$\ast$ given respectively by
$$
f\ast g(t)=\int_0^tf(t-s)g(s)ds, \quad t\geq0, \qquad
f\ast g(t)=\int_{-\infty}^{\infty} f(t-s)g(s)ds, \qquad t\in \R.
$$
These algebras  are canonical to define some algebra
homomorphisms (defined by integral representations) into
$\mathcal{B}(X).$ See further details in \cite{Ga-Mi-06}. The parameter $\alpha\ge 0$ denotes the fractional order of derivation of the Lebesque space; for $\alpha \in \NN$ is the iterative usual derivation and for $\alpha = 0$, the usual $p$-Lebesgue spaces.

Let's focus back on the generalized Stieltjes operators $\Sbmu$. Since the above described Sobolev-Lebesgue spaces are (continuously) contained in $L^p$, a natural question arises: are these operator $\mathcal{S}_{\beta, \mu}$ actually bounded in these spaces $\Tpa$ and $\Tpabs$? And in that case, which properties may they have as bounded operators? These questions depict the main motivation of this paper, which therefore is the study of the generalized Stieltjes operators $\mathcal{S}_{\beta, \mu}$ acting on $\mathcal{T}_p^{(\alpha)}(t^\alpha)$ and $\mathcal{T}_p^{(\alpha)}(|t|^\alpha)$. The powerful theory of subordination to $C_0$-groups of operators given in (\ref{subor}) allows us to get results regarding the boundedness, show commutative and factorization properties with the generalized Ces\`aro operator or to describe the spectrum set via a spectral mapping theorem. Also, we present some results involving the Hilbert and the Fourier transform, a new H$\ddot{\hbox{o}}$lder inequality on $\Tpa$, and a convolution-type theorem.

\section*{Outline and main results}

To present these results of the family of operators $\mathcal{S}_{\beta, \mu}$, the outline of this paper has been set as follows. In Section 2, two family of exponential functions $\phi_{\beta,\mu}$ and $\psi_{\gamma, \nu}$ which belong to $L^1(\R)$ for properly chosen parameters, are considered. Some of their properties are derived, especially those regarding their norms on $L^p(\R)$ spaces, its Fourier transform (some Beta functions) as well as the calculation of some useful convolutions, in particular, for $\beta,\gamma>0$ and $\lambda>\beta$, we have that
$$
\phi_{\beta, \lambda+\gamma}\ast\psi_{\gamma, \lambda-\beta}=\gamma B(\gamma,\lambda)\phi_{\beta,\lambda},
$$
(Proposition \ref{convolu}). These functions will play a main role to subordinate both generalized Stieltjes and Ces\`{a}ro operators in terms of $C_0$-groups.

In Section 3, we revisit some basic properties of the Sobolev-Lebesgue $\mathcal{T}_p^{(\alpha)}(t^\alpha)$ for $\alpha \ge 0$ and $p\ge 1$ included in some previous papers, mainly in \cite{[LMPS]}. We prove  a new  H\"{o}lder inequality in
Theorem \ref{Holder}:   $fg\in \mathcal{T}_1^{(\alpha)}(t^\alpha)$  with $f\in \mathcal{T}_p^{(\alpha)}(t^\alpha)$, $g\in \mathcal{T}_{p'}^{(\alpha)}(t^\alpha)$ and  $p, p'\in (1, \infty)$ conjugate exponents.

The generalized Stieltjes operator $\mathcal{S}_{\beta, \mu}$ acting on the Sobolev-Lebesgue $\mathcal{T}_p^{(\alpha)}(t^\alpha)$ is analyzed in section 4. First of all, we are able to subordinate $\Sbmu$ in terms of the $C_0$-group $(T_{t,p})_{t \in \R}$, via the exponential functions $\phi_{\beta,\mu}$ (Section 2), as an operator on $\Tpa$, that is
$$
\mathcal{S}_{\beta,\mu} f=\displaystyle\int_{-\infty}^\infty
\phi_{\mu -\beta +1/p, \mu}(r)T_{r,p}fdr=\theta(\phi_{\mu -\beta +1/p, \mu})f, \quad f\in \mathcal{T}_p^{(\alpha)}(t^\alpha).
$$

This immediately shows that $\Sbmu$ is bounded on $\Tpa$ and
$\|\mathcal{S}_{\beta, \mu}\| =  B(\beta-1/p,\mu-\beta+1/p)$
whenever $\alpha \geq 0$ and $0<\beta -1/p < \mu$ (Theorem \ref{cotaTa}). Its spectrum set, as operators on  $\mathcal{T}_p^{(\alpha)}(t^\alpha)$ for $p\ge 1, \: \alpha \geq 0$, is  obtained in Theorem \ref{spec},
$$
\sigma(\mathcal{S}_{\beta,\mu})= {\left\{B(\mu - \beta +1/p-i\xi, \beta -1/p+i\xi): \quad \xi \in \R\right\}}\cup\{0\}.
$$
This result extends the following ones, both cases stated as operators on $L^2(\R^+)$. First, $\sigma({\mathcal S})=[0,\pi]$, which was originally proved by Carleman in \cite[p. 169]{Ca}, and second, $\sigma(\mathcal{S}_{\beta,2\beta-1})=[0,B(\beta -1/2, \beta -1/2)]$ for $\beta>{1\over 2}$ in \cite[Proposition 1.1]{Federer}.

We are also able to give the adjoint of $\mathcal{S}_{\beta,\mu}$, which is precisely $\mathcal{S}_{\mu-\beta+1,\mu} \in \mathcal{B}(\mathcal{T}_{p'}^{(\alpha)}(t^{\alpha}))$, where $1/p + 1/p' = 1$, obtaining an elegant relation of the generalized Stieltjes operator with its adjoint (Theorem \ref{dualStieltjes}). Finally, we give an explicit expression for the composition $\mathcal{S}_{\beta,\mu} \mathcal{S}_{\gamma,\nu}$ involving the hypergeometric Gaussian function $_2F_1$, see details in Theorem \ref{compConDual2}.

In section 5, we study some properties of the generalized Stieltjes operator $\Sbmu$ regarding other operators, specifically the generalized Ces\'aro operator and the Hilbert transform on $\Tpa$. First, note that in \cite[Theorem 3.3]{[LMPS]}, one can write the generalized Ces\`aro operator as $\mathcal{C}_{\gamma} = \theta(\psi_{\gamma,1-1/p})$. Then we deduce straightforwardly that
 $\mathcal{S}_{\beta,\mu} \mathcal{C}_{\gamma} =
\mathcal{C}_{\gamma} \mathcal{S}_{\beta,\mu}$. We  are also able to find an explicit expression of $\mathcal{S}_{\beta, \mu} \mathcal{C}_\gamma$ in terms again of the hypergeometric Gaussian function $_2F_1$. Moreover, if $\beta = \gamma + 1$, the following elegant factorization holds
\begin{eqnarray*}
\mathcal{S}_{\gamma +1, \mu} \mathcal{C}_{\gamma} = \gamma B(\gamma, \mu - \gamma) \mathcal{S}_{1, \mu - \gamma}, \qquad \gamma >0, \,\quad \mu>\gamma+1-{1\over p},
\end{eqnarray*}
which it seems to be new (Theorem \ref{CesStieltComp}). The (semifinite) Hilbert transform ${\mathcal H}_+$, given by
$$
{\mathcal H}_+f(t):=\hbox{v.p.}{i\over \pi}\int_0^\infty{f(s)\over t-s}ds, \qquad t>0,\qquad f\in \mathcal{T}_p^{(\alpha)}( t^\alpha),
$$
 has been studied as a bounded operator on $L^p(\R^+)$ (see for example \cite[Section 6.1]{[HKV]}). Then, we prove that it is also a bounded operator on $\mathcal{T}_p^{(\alpha)}( t^\alpha)$ which commutes with the generalized Stieltjes transform $\mathcal{S}_{\beta,\mu}$ for $1<p<\infty$ and $\alpha \ge 0$, see Theorem \ref{semihilbert}.

In Section 6, we introduce generalized Stieltjes operators ${\mathcal S}_{\beta, \mu}$ on the
Sobolev-Lebesgue spaces $\mathcal{T}_p^{(\alpha)}(\vert t\vert^\alpha)$ on $\mathbb{R}.$
Again, it is relevant to mention that the $C_0$-group of
isometries on $\mathcal{T}_p^{(\alpha)}(\vert t\vert^\alpha)$,
$(T_{t, p})_{t\in \RR}$
 is the main tool to prove the some of main results
in this section (Theorem \ref{lemma22.1}).

Note that Stieltjes operator ${\mathcal S}$ is closer to Hilbert transform ${\mathcal H}$, $
{\mathcal H}f(-t)= -{\mathcal S}f(t),$ for $t>0$, a.e. where
$$
{\mathcal H}f(t):=\hbox{v.p.}{i\over \pi}\int_{-\infty}^\infty {f(s)\over t-s}ds, \qquad t\in \R, \,\,a.e,
$$
is a bounded operator from $L^p(\R)$ onto $L^p(\R)$ for $1<p<\infty$ and an isometry for $p=2$, see a real proof in \cite{[Du]}. We prove that this also holds in $\Tpabs$ for $\alpha \geq 0$, that is, the operator ${\mathcal H}$ is bounded on  $\mathcal{T}_p^{(\alpha)}(\vert t\vert^\alpha)$ for $1<p<\infty$, an isometry for $p=2$
(Theorem \ref{Hilberttransf}), and that it commutes with $\Sbmu$. Finally the connection with Fourier transform  is presented,
$$
\widehat{{\mathcal S}_{\beta, \mu}(f)}={\mathcal S}_{\mu-\beta+1,\mu}(\widehat{f}), \qquad
f\in L^p(\R),
$$
for $1\le p\le 2$ and  $0< \beta -1/p < \mu$, (Theorem \ref{conmutan2}).

In section 7, we take a look at the convolution $\otimes$ defined from $L^p(\R^+)\times L^q(\R^+)$ to $L^r(\R^+), \: \: 1/p +1/q = 1/r$,
$$
f \otimes g  :=f{\mathcal H}_+g + g{\mathcal H}_+f, \qquad f\in L^p(\R^+), \,g\in L^q(\R^+),
$$
which can be found in \cite{[SV],[Ya],[YaM]}, and where  the following interesting relation involving the Stieltjes transform $\mathcal{S}(f\otimes g) = (\mathcal{S}f)(\mathcal{S}g)$ is proven. We wonder if  the generalized Stieltjes operator $\mathcal{S}_{\beta, \mu}$ also behaves in such elegant way. In fact, we are able to obtain a nice-looking formula for $\mathcal{S}_{n,m}(f\otimes g)$, involving a sum of products of different generalized Stieltjes operators, whenever $1\leq n \leq m$ and both $n,m$ are natural numbers (Theorem \ref{convProd}).

Finally in the last section, we use the software Mathematica  to visualize the spectrum  set of the generalized Stieltjes operators, $\sigma(\mathcal{S}_{\beta,\mu}),$ on ${\mathcal T}_p^{(\alpha)}(t^\alpha)$ in some particular cases. We also present some ideas to continue this research in subspaces of $H^p(\C^+)$ for $1<p<\infty$ or some problems which remain open after this research.

\bigskip

\noindent{\bf Notation.} For $1\le p<\infty$  recall that the Banach space $L^p(\mathbb{R}^+)$ is the set of
Lebesgue $p$-integrable functions, that is, $f$ is a measurable
function and
$$
||f||_p:=\left(\int_0^\infty |f(t)|^pdt\right)^{1/p}<\infty.
$$
The space $L^\infty(\mathbb{R}^+)$ be the set of  measurable functions such that
$$
||f||_\infty:=\hbox{ess sup}\{\vert f(t)\vert \,\,:\,\,t\in \R\}<\infty.
$$
In the case that $f$ is a continuous function, then $||f||_\infty=\hbox{sup}\{\vert f(t)\vert \,\,:\,\,t\in \R\}.$

 Recall that the Gamma and Beta functions, also called the Euler
integrals of the first kind, are defined by:
$$
\Gamma(z)= \int_0^\infty t^{z-1} e^{-t} dt,\qquad B(z,w)= \int_0^1 t^{z-1} (1-t)^{w-1} dt, \quad \Re z,\, \Re w>0,
$$
and satisfy the property $ B(z,w)= \frac{\Gamma(z)
\Gamma(w)}{\Gamma(z+w)}$, where $\Gamma$ denotes the Gamma function.

\section{Two families  of exponential funtions on $L^1(\R)$}

Let introduce two families of parametrized functions, which will be denoted by $\phi_{\beta,\mu}$ and $\psi_{\gamma,\nu}$. Our interest in these functions is due to that they appear in the integral representations of the generalized Stieltjes operator $\Sbmu$ and generalized Ces\`aro operator $\mathcal{\gamma}$ respectively, as we shall se in Sections 4 and 5.

We define the set of functions $(\phi_{\beta,\mu})_{\beta,\mu \in \R}$ by
\begin{equation}\label{functionp}
\phi_{\beta, \mu}(t):={e^{\beta t}\over (1+e^{t})^\mu}, \qquad t\in \R.
\end{equation}
Note that $
\phi_{\beta, \mu}(-t)=
\phi_{\mu-\beta, \mu}(t)
$ for $t\in \R$ and $\phi_{\beta, \mu}\phi_{\gamma, \nu}=\phi_{\beta+\gamma, \mu+\nu}$ for $\beta, \gamma,\mu,\nu\in \R$. It is direct to check that
\begin{eqnarray*}
\phi'_{\beta, \mu}&=& \phi_{\beta, \mu}(\beta-\mu \phi_{1,1}),\cr
\phi''_{\beta, \mu}&=& \phi_{\beta, \mu}(\beta^2-\mu(2\beta +1) \phi_{1,1}+\mu(\mu+1)\phi^2_{1,1}),
\end{eqnarray*}
and, in general, $\phi^{(n)}_{\beta, \mu}= \phi_{\beta, \mu}P_{n,\beta,\mu}(\phi_{1,1})$ where $P_{n,\beta,\mu}$  is a polynomial of degree $n$ whose coefficients depends on $\beta$ and $\mu$, that is, $P_{n,\beta,\mu} (z) = \sum_{m=0}^{n}{a_{m}^{(n,\beta,\mu)} z^{m}}$ with
\begin{eqnarray*}
a_{m}^{(n,\beta,\mu)} = {{\mu+m-1}\choose m}\sum_{j=0}^{m}{(-1)^{j}{m\choose j}(\beta +j)^{n}},\qquad	n\geq m\geq 0.
\end{eqnarray*} To show this, note that $\phi_{\beta,\mu}P_{n+1}(\phi_{1,1}) = \phi_{\beta,\mu}^{(n+1)} = ((\phi_{\beta,\mu}P_{n}(\phi_{1,1}))'$, and one can derive the subsequent recursive formula:
\begin{eqnarray*}
a_{m}^{(n+1,\beta,\mu)}-(\beta +m)a_{m}^{(n,\beta,\mu)}=-(\mu +m-1)a_{m-1}^{(n,\beta,\mu)}, \qquad  n \geq m \geq 0,
\end{eqnarray*}
and $a_{0}^{(0,\beta,\mu)}=1$, $a_{-1}^{(n,\beta,\mu)}=0$, $a_{n+1}^{(n,\beta,\mu)}=0,$ for  $n \geq 0$. Then, one can show by induction on $m$ that the given closed-form expression of coefficients $a_{m}^{(n,\beta,\mu)}$ is the (unique) solution of this recursive formula.

The Gaussian hypergeometric function will be needed for our next result involving the family of functions $(\phi_{\beta,\mu})_{\beta, \mu \in \R}$. As usual, we denote it by $_2F_1$, and it is given by
$$
_2F_1(a,b;c;z):={\Gamma(c)\over
	\Gamma(b)\Gamma(a)}\sum_{n=0}^\infty{\Gamma(a+n)\Gamma(b+n)\over
	\Gamma(c+n)}{z^n\over n!}, \qquad \vert z\vert<1,
$$
(\cite[Section 9.1]{[GR]}). A special case occurs when $a=c$, then:
$_2F_1(c,b;c;z) = (1-z)^{-b}$. Even more, for $\Re c>\Re b>0,$ it can
be analytically  extended  via the  integral
\begin{eqnarray}\label{extendedGaussian}
	_2F_1(a,b;c;z) :=\frac{\Gamma(c)}{\Gamma(c-b)\Gamma(b)} \int_0^1 s^{b-1}(1-s)^{c-b-1}(1-zs)^{-a}ds, \qquad z\in \mathbb{C} \backslash [1,+\infty),
\end{eqnarray}
(\cite[Formula 9.111]{[GR]}).
We will also make use of the following identity:
\begin{eqnarray}\label{kummer}
_2F_1(a,b;c;z) = (1-z)^{-a} {}_2F_1\left(a,c-b;c;{z\over{z-1}}\right)
\end{eqnarray}
(\cite[Formula 9.131 (1)]{[GR]}).

\begin{theorem} \label{cottas}Fixed $\beta, \mu\in \R$.
\begin{itemize}
\item[(i)] $\phi_{\beta, \mu }\in C_b(\R)$ if and only if  $0\le \beta\le \mu$ and $\Vert \phi_{0, 0}\Vert_\infty=\Vert \phi_{0, \mu}\Vert_\infty=\Vert \phi_{\beta, \beta}\Vert_\infty=1$,
$$\Vert \phi_{\beta, \mu}\Vert_\infty= \left({\mu-\beta\over \mu}\right)^\mu\left({\beta\over \mu-\beta}\right)^\beta,$$
for $0<\beta<\mu$.
\item[(ii)] For $p\ge 1$, $\phi_{\beta, \mu}\in L^p(\R)$ if and only if  $0< \beta< \mu$  and
$$
\Vert \phi_{\beta, \mu}\Vert_p=\left(B(p\beta, p(\mu-\beta))\right)^{1\over p}.
$$
\item[(iii)]  For  $0< \beta< \mu$ , we have that
$$
\widehat{\phi_{\beta, \mu}}(\xi)=B(\beta-i\xi, \mu-\beta+i\xi), \qquad \xi \in \R.
$$
\item[(iv)]Let $\beta, \mu, \gamma, \nu \in \R$ be such that $0<\beta < \mu, 0<\gamma < \nu$. Then
\begin{eqnarray*}
\phi_{\beta,\mu} * \phi_{\gamma, \nu} (t) = B(\mu -\beta + \gamma, \nu -\gamma + \beta)
e^{-(\mu -\beta)t} {}_2F_1(\mu, \mu -\beta + \gamma ; \mu + \nu; 1-e^{-t}),
\end{eqnarray*}
for $ t\in \R.$
\end{itemize}
\end{theorem}

\bgproof (i)
It is straightforward to conclude that  $\phi_{\beta, \mu }\in C_b(\R)$
if and only if  $0\le \beta\le \mu$. In this case,
we have that $0< \phi_{\beta, \mu }(t) \leq 1,$ for  $t \in \R$. Therefore
$\Vert \phi_{0, 0}\Vert_\infty=\Vert \phi_{0, \mu}\Vert_\infty=\Vert \phi_{\beta, \beta}\Vert_\infty=1$.
If $0 < \beta < \mu$, then $\phi_{\beta, \mu}'(t) = 0 $ if and only if $ e^t = {\frac{\beta}{\mu-\beta}}$,
which is the point where $\phi_{\beta,\mu}$
reaches its maximum value, $\|\phi_{\beta,\mu}\|_{\infty}$, claimed in the statement.

(ii) By (i), it is clear that $\phi_{\beta, \mu}\in L^p(\R)$
implies that  $0< \beta< \mu$. Under this latter assumption, straightforward calculations show that actually
$\phi_{\beta, \mu}\in L^p(\R$) and
$$
\|\phi_{\beta, \mu}\|_p = \left(\int_0^{\infty} \frac{s^{\beta p -1}}{(1+s)^{\mu p}}ds \right)^{1\over p}
= (B(p\beta, p(\mu-\beta)))^{1\over p}.
$$
(iii) For $\xi \in \R$, we have that
\begin{eqnarray*}
\widehat{\phi_{\beta, \mu}}(\xi) = \int_{0}^{\infty}\frac{s^{\beta -i\xi -1}}{(1+s)^{\mu}}ds
=B(\beta-i\xi,\mu-\beta+i\xi),\qquad \xi\in \R.
\end{eqnarray*}
(iv) Let $0<\beta < \mu$, and $0<\gamma < \nu$. Then, by \cite[p. 317, 3.197(1)]{[GR]} and \ref{kummer}, we have that
\begin{eqnarray*}
\phi_{\beta,\mu} * \phi_{\gamma, \nu} (t) &=& \int_{-\infty}^{\infty}
\frac{e^{\beta \tau}}{(1+e^\tau)^\mu} \frac{e^{\gamma(t-\tau)}}{(1+e^{t-\tau})^\nu} d\tau
= e^{\gamma t} \int_0^{\infty} \frac{s^{\nu - \gamma + \beta -1}}{(1+s)^\mu (e^t +s)^\nu} ds
\\
&=& e^{\beta t}  B(\mu -\beta + \gamma, \nu -\gamma +\beta) {}_2F_1(\mu, \nu -\gamma +\beta ;\mu + \nu; 1-e^{t})
\\
&=& B(\mu -\beta + \gamma, \nu -\gamma +\beta)e^{-(\mu -\beta)t} {}_2F_1(\mu, \mu -\beta +\gamma ;\mu + \nu; 1-e^{-t}),
\end{eqnarray*}
for $t\in \R$ and we conclude the proof.
\edproof

Some particular cases of $\phi_{\beta,\mu} * \phi_{\gamma, \nu}$ are given in the next result. The equality in part (ii) is given in \cite[p. 317, 3.197(7)]{[GR]}.
\begin{corollary}
\begin{itemize}
\item[(i)] For $0<\beta <\mu$, we have
$$
\phi_{\beta,\mu}\ast \phi_{\beta,\mu}(t)= B(\mu, \mu)e^{-(\mu-\beta)t}{}_2F_1(\mu, \mu  ;2\mu; 1-e^{-t}),\qquad t\in \R,
$$
in particular $\phi_{\beta,1}\ast \phi_{\beta,1}(t)=\displaystyle{t e^{-(1-\beta)t}\over 1-e^{-t}}$ for $t\in \R\backslash\{0\}$ and $0<\beta<1$.
\item[(ii)]
Let $0<\mu$. Then
$$\phi_{\beta+{1\over 2}, \mu} \ast \phi_{\beta,\mu}(t) = \sqrt{\pi} { \Gamma(\mu -{1 \over 2})\over \Gamma(\mu)}\phi_{2\beta, 2\mu-1}\left({t\over 2}\right), \qquad t\in \R.
$$
\end{itemize}
\end{corollary}

For $\gamma, \nu \in \R$, we consider a second family of functions,
\begin{eqnarray} \label{psifunctions}
	\psi_{\gamma, \nu}(t):= \gamma (1-e^{-t})^{\gamma-1} e^{-\nu t}\chi_{(0,\infty)}(t), \qquad t\in \R,
\end{eqnarray}
where $\chi_{(0,\infty)}$ is the characteristic funtion on the interval $(0,\infty)$. Note $\psi_{1, \nu}=e_{\nu}$ where $e_{\nu}(t):=e^{-\nu t} \chi_{0,\infty}(t)$ for $\nu,t \in \R$. In the last part of this section, we give first an analogous result to Theorem \ref{cottas} for the family of functions $(\psi_{\gamma, \nu})_{\gamma, \nu \in \R}$. After that, we provide explicitly some useful convolutions involving the two families $(\phi_{\beta,\mu})_{\beta, \mu \in \R}$ and $(\psi_{\gamma, \nu})_{\gamma, \nu \in \R}$.

\begin{proposition} Fixed $\gamma, \nu \in \R$.
\begin{itemize}
\item[(i)] $\psi_{\gamma, \nu }\in C_b(\R)$ if and only if  $1 \leq \gamma$,  $0 \leq \nu$. In that case, $\|\psi_{1, \nu}\|_{\infty} = 1,  \|\psi_{\gamma, 0}\|_{\infty} = \gamma$, and
$$\Vert \psi_{\gamma, \nu}\Vert_\infty= \gamma \left(\frac{\gamma -1}{\nu +\gamma -1}\right)^{\gamma-1}
\left(\frac{\nu}{\nu +\gamma -1}\right)^{\nu}$$
if $1 < \gamma$,  $0 < \nu$.
\item[(ii)] For $p\ge 1$, $\psi_{\gamma, \nu}\in L^p(\R)$ if and only if $0< \nu$ and also, either $\gamma\ge 1$, or $p<{1\over 1-\gamma}$ for $0<\gamma<1$; in both cases
$$
\Vert \psi_{\gamma, \nu}\Vert_p=\gamma \left(B(p\nu, p(\gamma-1)+1)\right)^{1\over p}.
$$
\item[(iii)]  For  $\gamma,\nu>0$, we have that
$$
\widehat{\psi_{\gamma, \nu}}(\xi)=\gamma B(\gamma, \nu+i\xi), \qquad \xi \in \R.
$$
\end{itemize}
\end{proposition}

\bgproof (i) It is clear that $\psi_{\gamma, \nu }\in C_b(\R)$ if and only if  $1 \leq \gamma$,  $0 \leq \nu$.
Moreover, since under these assumptions $0 \leq \psi_{\gamma, \nu } (t) \leq \gamma$, we have that
$\|\psi_{1, \nu}\|_{\infty} = 1,  \|\psi_{\gamma, 0}\|_{\infty} = \gamma$.
If we consider $1 < \gamma, 0 < \nu$, and $t>0$, then $\psi_{\gamma, \nu }'(t) = 0$ if and only if
$e^{-t} = \frac{\nu}{\nu + \gamma -1}$, where it reaches the maximum $\|\psi_{\gamma, \nu }\|_{\infty},$
given in the statement.

(ii)  For $p \geq 1$, we apply the  change of variable $s = e^{-t}$ to get
\begin{eqnarray*}
\Vert \psi_{\gamma, \nu}\Vert_p =
\gamma \left(\int_0^1(1-s)^{p(\gamma -1)}s^{p\nu -1}ds\right)^{1 \over p}=
\gamma \left(B(p\nu, p\gamma-p+1)\right)^{1\over p},
\end{eqnarray*}
whenever $\nu >0$ and $p\gamma - p +1 >0$, that is, either $\gamma \geq 1$ or $0< \gamma <1$ and
$p < \frac{1}{1-\gamma}$.

(iii) Applying the change of variable $s = e^{-t}$  and assuming $\gamma,\nu>0$:
\begin{eqnarray*}
\widehat{\psi_{\gamma, \nu}}(\xi) = \int_0^1 \gamma (1-s)^{\gamma -1} s^{\nu +i\xi -1}ds
= \gamma B(\gamma, \nu+i\xi),
\end{eqnarray*}
for $\xi \in \R$ and we conclude the proof.
\edproof

\begin{proposition}\label{convolu} Fixed $\beta,\gamma, \mu, \nu >0$. Then
$$
\phi_{\beta, \mu}\ast\psi_{\gamma, \nu}(t)=\gamma B(\gamma,\beta+\nu)\phi_{\beta,\mu}(t)_2F_1\left(\mu,\gamma;\beta+ \gamma +\nu;{e^t\over e^t+1}\right), \qquad t\in \R.
$$
In particular,
\begin{itemize}
\item[(i)] for $\beta,\gamma>0$ and $\lambda>\beta$, we have that
\begin{equation}\label{factor}
\phi_{\beta, \lambda+\gamma}\ast\psi_{\gamma, \lambda-\beta}=\gamma B(\gamma,\lambda)\phi_{\beta,\lambda}.
\end{equation}
\item[(ii)] for $0<\beta <1$ and $t\in \R$,
\begin{eqnarray*}
\phi_{\beta, \mu}*\psi_{1,1-\beta}(t) &=& {{e^{-(1-\beta)t}}\over{\mu -1}} \left(1-\frac{1}{(1+e^t)^{\mu-1}}\right), \qquad \mu \ne 1 \\
\phi_{\beta, 1}*\psi_{1,1-\beta} (t)&=& e^{-(1 -\beta)t} \log(1+e^t).
\end{eqnarray*}
\end{itemize}
\end{proposition}

\bgproof
 We apply the change of variable $e^{-\tau} = s$, and by identities (\ref{extendedGaussian}) and (\ref{kummer}), we get that
 \begin{eqnarray*}
\phi_{\beta, \mu} * \psi_{\gamma ,\nu}(t) &=& \gamma  e^{\beta t} \int_0^1 s^{\nu +\beta -1} (1-s)^{\gamma -1}(1+e^{t}s)^{-\mu}ds
\\
&=& \gamma  e^{\beta t} B(\gamma, \beta + \nu)
{}_2 F_1(\mu ,\beta + \nu ; \beta+ \gamma +\nu; -e^t)
\\
&=& \gamma  B(\gamma, \beta + \nu)
 \frac{e^{\beta t}}{(1+e^t)^{\mu}} {}_2 F_1\left(\mu ,\gamma ; \beta+ \gamma +\nu; \frac{e^t}{1+e^t}\right)
\\
&=& \gamma B(\gamma, \beta+ \nu)\phi_{\beta,\mu}(t)
 {}_2 F_1\left(\mu ,\gamma ; \beta+ \gamma +\nu; \frac{e^t}{1+e^t}\right),
\end{eqnarray*}
for any $t \in \R$. Since, as pointed out at the beginning of this section, $_2F_1(c,b;c;z) = (1-z)^{-b}$,  we have that:
$$
\phi_{\beta, \lambda+\gamma}\ast\psi_{\gamma, \lambda-\beta} (t)
 = \gamma B(\gamma,\lambda)\phi_{\beta,\lambda + \gamma}(t) \left(1 - \frac{e^t}{1+e^t}\right)^{-\gamma}
 = \gamma B(\gamma,\lambda)\phi_{\beta,\lambda}(t)
$$
for any $t \in \R$ and we show (i). Both equalities in (ii) can be directly obtained as particular cases of the integral at the beginning of this proof.  \edproof

\section{Sobolev-Lebesgue spaces $\mathcal{T}_p^{(\alpha)}(t^\alpha)$ }

\setcounter{theorem}{0} \setcounter{equation}{0}

Let $\mathcal{D}_+$ be the class of $C^\infty$-functions with
compact support on $[0,\infty)$ and $\mathcal{S}_+$ the Schwartz
class on $[0,\infty)$. For a function $f\in \mathcal{S}_+$ and
$\alpha>0$, the \textit{Weyl fractional integral} of order
$\alpha$, $W^{-\alpha}_+f$, is defined by
$$
W^{-\alpha}_+f(t):=\frac{1}{\Gamma(\alpha)}\int_t^\infty
(s-t)^{\alpha-1}f(s)ds, \qquad t\in\mathbb{R}^+.
$$
The \textit{Weyl fractional derivative} $W^\alpha_+f$ of order
$\alpha$ is defined by
$$
W^\alpha_+ f(t):=(-1)^n\frac{d^n}{dt^n}W_+^{-(n-\alpha)}f(t),\quad
t\in \mathbb{R}^+
$$
where $n=[\alpha]+1,$ and $[\alpha]$ denotes the integer part of
$\alpha.$ It is proved that
$W_+^{\alpha+\beta}=W_+^{\alpha}(W_+^{\beta})$ for any $\alpha,
\beta\in\mathbb{R},$ where $W_+^0=Id$ is the identity operator, and
$(-1)^nW_+^n=\frac{d^n}{dt^n}$ holds with $n\in\mathbb{N}$, see
more details in \cite{Mi-Ro, Sa-Ki-Ma-93}.

%Given $f\in C(\mathbb{R};X),$  the Riemann-Liouville integral  $D^{-\alpha}f,$ of order $\alpha > 0$ is defined by
%$$
% D^{-\alpha}f(t):=\displaystyle{\small\int_0^t \frac{(t-s)^{\alpha-1}}{\Gamma(\alpha)}}f(s)ds,\,\mbox{ if } t\geq 0.
%$$

Take $ \lambda>0$ and define $f_\lambda$ by
$f_\lambda(r):=f(\lambda r)$ for $r>0$ and $f \in \mathcal{S}_+$.
It is direct to check that
\begin{equation}\label{escala}
W^\alpha_+f_\lambda = \lambda^\alpha (W^\alpha_+f)_\lambda, \qquad
f \in \mathcal{S}_+,
\end{equation}
for $\alpha \in \RR$.

%We also recall, the Riemann-Liouville fractional calculus. Let $\alpha>0$ and $f$ continuos on $\mathbb{R}$. The \textit{Riemann-Liouville fractional integral}, $D_0^\alpha f$ of order $\alpha$, is defined by

%$$
%D_0^{-\alpha} f(t)= \left\{ \begin{array}{l}
% D^{-\alpha}_+f(t):=\displaystyle\frac{1}{\Gamma(\alpha)}\int_0^t (t-s)^{\alpha-1}f(s)ds,\,\mbox{ if } t\geq 0, \\
% D^{-\alpha}_-f(t):=\displaystyle\frac{1}{\Gamma(\alpha)}\int_t^0 (s-t)^{\alpha-1}f(s)ds,\,\mbox{ if } t< 0. \\
% \end{array} \right.
%$$

%For $n\in \mathbb{N}\cup \{0\}\cup\{\infty\}$, let $C^{(n)}(\mathbb{R})$ be the set of all functions $n$-times differentiable on $\mathbb{R}$. Let $n=[\alpha]+1$, where $[\alpha]$ denote the integer part of $\alpha$. If $f\in C^{(n)}(\mathbb{R})$, we define the Riemann-Liouville  fractional derivative of order $\alpha>0$ by

%$$ D_0^\alpha f(t)= \left\{ \begin{array}{l}
% D^{\alpha}_+f(t):=\displaystyle(-1)^n\frac{d^n}{dt^n}D_+^{-(n-\alpha)}f(t),\,\mbox{ if } t\geq 0, \\
%  D^{\alpha}_-f(t):=\displaystyle(-1)^n\frac{d^n}{dt^n}D_-^{-(n-\alpha)}f(t),\,\mbox{ if } t< 0. \\
% \end{array} \right.
%$$

Now we recall a family of subspaces
$\mathcal{T}_p^{(\alpha)}(t^\alpha)$ which are contained in
$L^p(\RR^+)$ with $1\le p < \infty$.
For $\alpha>0$,  the Banach space
$\mathcal{T}_p^{(\alpha)}(t^\alpha)$ is defined as the completion of
the Schwartz class  $\mathcal{S}_+$ in the norm
\begin{equation*}
||f||_{\alpha, p}:=\frac{1}{\Gamma(\alpha+1)}\left(\int_0^\infty
|W^\alpha_+ f(t)|^pt^{\alpha p}dt\right)^{\frac{1}{p}},
\end{equation*}
 \cite[Definition 2.1]{[LMPS]}. As usual, We write  $\mathcal{T}_p^{(0)}(t^0)= L^p(\mathbb{R}^+)$
and $||\quad||_{0, p}=||\quad||_{p}$.  The case $p=1$ and $\alpha
\in \NN$ where introduced in \cite{AK} and for $\alpha >0$ in
\cite{Ga-Mi-06}.

These spaces $ \mathcal{T}_p^{(\alpha)}(t^\alpha)$ have similar properties than the Lebesgue spaces $L^p(\RR^+)$,  see \cite{Ro-08} and \cite[Proposition 2.2]{[LMPS]}. In fact  the operator $D^\alpha_+: \mathcal{T}_p^{(\alpha)}(t^\alpha)\to L^p(\RR^+)$  defined by
$$
f\mapsto D^\alpha_+ f(t) = \frac{1}{\Gamma(\alpha+1)} t^\alpha
W^\alpha_+f(t), \qquad t\ge 0, \quad f \in
\mathcal{T}_p^{(\alpha)}(t^\alpha).
$$
is an isometry, i.e., $\Vert f\Vert_{\alpha,p}=\Vert D^\alpha_+ f\Vert_p$ and the inverse operator $(D^\alpha_+)^{-1}: L^p(\RR^+)\to \mathcal{T}_p^{(\alpha)}(t^\alpha)$ is given by
\begin{equation}\label{inverse}
(D^\alpha_+)^{-1} f(t):={\alpha}\int_t^\infty (s-t)^{\alpha-1}{f(s)\over s^\alpha} ds, \qquad t>0,\qquad f\in  L^p(\RR^+).
\end{equation}

 For  $p\geq 1$ and $\beta>\alpha>0$, we also have that
\begin{itemize}
\item[(i)] $\mathcal{T}_p^{(\beta)}(t^\beta)\hookrightarrow \mathcal{T}_p^{(\alpha)}(t^\alpha)\hookrightarrow  L^p(\mathbb{R}^+) $, with the inclusions being continuous maps.

\item[(ii)] $\mathcal{T}_p^{(\alpha)}(t^\alpha)\ast \mathcal{T}_1^{(\alpha)}(t^\alpha) \hookrightarrow \mathcal{T}_p^{(\alpha)}(t^\alpha)$ for $1\le p <\infty$, where
\begin{equation}\label{convos}
f\ast g(t)=\int_0^tf(t-s)g(s)ds, \quad t\ge 0, \qquad f\in
\mathcal{T}_p^{(\alpha)}(t^\alpha), \quad
g\in\mathcal{T}_1^{(\alpha)}(t^\alpha).
\end{equation}

\item[(iii)] If $p>1$ and $p'$ satisfies $\frac{1}{p}+\frac{1}{p'}=1$, then the dual of $\mathcal{T}_p^{(\alpha)}(t^\alpha)$ is $\mathcal{T}_{p'}^{(\alpha)}(t^\alpha)$, where the duality is given by $$\langle f,g \rangle_{\alpha} ={1\over \Gamma(\alpha+1)^2}\int_0^\infty W^\alpha_+ f(t)W^\alpha_+ g(t)t^{2\alpha}dt=\langle D^\alpha_+ f,D^\alpha_+ g
\rangle_{0},$$
for $f\in \mathcal{T}_p^{(\alpha)}(t^\alpha)$, $g\in \mathcal{T}_{p'}^{(\alpha)}(t^\alpha)$.
\end{itemize}

If $\alpha,a>0$ and $p\geq 1,$ then

(i) $t^{\beta}\not \in \mathcal{T}_p^{(\alpha)}(t^\alpha)$ for
$\beta\in\mathbb{C}.$

(ii) $(a+t)^{-\beta}\in \mathcal{T}_p^{(\alpha)}(t^\alpha)$ for
${\Re}\beta>1/p.$
(\cite[Lemma 2.3]{[LMPS]}).

\begin{remark}\label{infinito}{\rm For $p=\infty$, the nature of the space $\mathcal{T}_\infty^{(\alpha)}(t^\alpha)$  is different than $\mathcal{T}_p^{(\alpha)}(t^\alpha)$ for $p\ge1$. One might think to define

$$
||f||_{\alpha,\infty}:={1\over \Gamma(\alpha+1)}\hbox{ess sup}\{t^\alpha\vert W^{\alpha}_+ f(t)\vert \,\,:\,\,t\in \R^+\}=||D_+^\alpha f||_{\infty}<\infty.
$$
Consider $l(t):=\log(1+t)$ ($t\in \R^+)$; note that $||l||_{1,\infty}<\infty$ and $\Vert l\Vert_\infty=+\infty$. However in the case that $0<\alpha<\beta$, note that for $t>0$,
$$
\vert D_+^\alpha f(t)\vert\le {t^{\alpha}\over \Gamma(\alpha+1)}\int_t^\infty{(s-t)^{\beta-\alpha-1}\over s^\beta\Gamma(\beta-\alpha)}s^\beta\vert W^{\beta}_+f(s)\vert ds\le {\beta\over \alpha}||f||_{\beta,\infty},
$$
and we conclude $||f||_{\alpha,\infty}\le  {\beta\over \alpha}||f||_{\beta,\infty}$ for $f$ in the Schwarz class on $\R^+$. Due to, we will skip the case $p=\infty$ on this paper, although some of next results are also valid in this particular case. }

\end{remark}
\medskip

In the next theorem, we extend the   H$\ddot{\hbox{o}}$lder inequality in Lebesgue spaces $L^p(\R^+)$ to the case of fractional Sobolev-Lebesgue spaces  $\mathcal{T}_p^{(\alpha)}(t^\alpha)$. Note that in the integer case $\alpha \in \NN$, it a straightforward consequence of Leibnitz formula. To attack the general case we use the following Leibnitz formula
\begin{eqnarray*}\label{Leib}
W^{\alpha}_+(fg)(x)=f(x)W^{\alpha}_+g(x)+g(x)W^{\alpha}_+f(x)-\int_x^\infty\int_x^\infty{d\varphi^{\alpha-1}_{t,u}(x)\over dx}W^{\alpha}_+f(t)W^{\alpha}_+g(u)dtdu
\end{eqnarray*}
for $f, g \in  {\mathcal S}_+$ and $\alpha>0$ \cite[Proposition 2.5]{[GP]}, where the function $\varphi^{\alpha-1}_{t,u}$ is given by
$$
\varphi^{\alpha-1}_{t,u} (x) = {(u-x)^{\alpha-1}\over \Gamma(\alpha)}\,_2F_1\left(-\alpha+1,-\alpha+1,1;{t-x\over u-x}\right),
$$
for $x<t<u$; $\varphi^{\alpha-1}_{t,u}(x)= \varphi^{\alpha-1}_{u,t}(x)$ for  $x<u<t$ and let $\varphi^{\alpha-1}_{t,u}=0$ in all other cases (\cite[p. 313]{[GP]}).

\begin{theorem}\label{Holder} (H$\ddot{\hbox{o}}$lder inequality)
Take $\alpha\ge 0,$ and $p, p'\in (1, \infty)$ conjugate exponents. Given $f\in \mathcal{T}_p^{(\alpha)}(t^\alpha)$ and $g\in \mathcal{T}_{p'}^{(\alpha)}(t^\alpha)$ then $fg\in \mathcal{T}_1^{(\alpha)}(t^\alpha)$ and
$$
\Vert fg\Vert_{\alpha,1}\le C_\alpha \Vert f\Vert_{\alpha,p}\Vert g\Vert_{\alpha,p'},
$$
where $C_\alpha$ is a positive constant.

\end{theorem}

\bgproof Take $f, g \in {\mathcal S}_+$. By the Leibniz formula, we have that
\begin{eqnarray*}
\Vert fg\Vert_{\alpha,1}\le{1\over \Gamma(\alpha +1)} \int_0^\infty x^\alpha \vert f(x)\vert \,\vert W^{\alpha}_+g(x)\vert dx+ {1\over \Gamma(\alpha +1)}\int_0^\infty x^\alpha \vert g(x)\vert \, \vert W^{\alpha}_+f(x)\vert dx\\
+{1\over \Gamma(\alpha +1)}\int_0^\infty x^\alpha \int_x^\infty\int_x^\infty\left|{d\varphi^{\alpha-1}_{t,u}(x)\over dx}\right|\vert W^{\alpha}_+f(t)\vert \,\vert W^{\alpha}_+g(u)\vert dtdudx.
\end{eqnarray*}

First, notice that by classical Hölder inequality, and the continuous inclusion $\Tpa \hookrightarrow L^p(\R^+)$, we have that
\begin{eqnarray} \label{Holderreason}
	{1\over \Gamma(\alpha +1)} \int_0^\infty x^\alpha \vert f(x)\vert \,\vert W^{\alpha}_+g(x)\vert dx
	\leq \|f\|_{0,p} \|g\|_{\alpha, p'} \leq C_\alpha \|f\|_{\alpha,p} \|g\|_{\alpha, p'}
\end{eqnarray}
For convenience and clarity of writing, we will denote all these constants by $C_\alpha$. The argument above let us bound the first two terms in the above sum by $C_\alpha \|f\|_{\alpha,p} \|g\|_{\alpha, p'}$, so it is sufficient that the inequality holds for the third one, which will be the aim for the rest of the proof. We apply twice Fubini theorem in the third addend, labelled as $A_\alpha(f,g)$, to get
\begin{eqnarray*}
A_\alpha(f,g)  &:=&\int_0^\infty x^\alpha \int_x^\infty\int_x^\infty\left|{d\varphi^{\alpha-1}_{t,u}(x)\over dx}\right|\vert W^{\alpha}_+f(t)\vert \,\vert W^{\alpha}_+g(u)\vert dtdudx\\
&\,&=\int_0^\infty \vert W^{\alpha}_+g(u)\vert  \int_0^u x^\alpha \int_x^\infty\left|{d\varphi^{\alpha-1}_{t,u}(x)\over dx}\right|\vert W^{\alpha}_+f(t)\vert \, dtdxdu\\
&\,&=\int_0^\infty \vert W^{\alpha}_+g(u)\vert  \int_0^u \vert W^{\alpha}_+f(t)\vert  \int_0^t x^\alpha \left|{d\varphi^{\alpha-1}_{t,u}(x)\over dx}\right| \, dxdtdu\\
&\,&\qquad +
\int_0^\infty \vert W^{\alpha}_+g(u)\vert  \int_u^\infty \vert W^{\alpha}_+f(t)\vert  \int_0^u x^\alpha \left|{d\varphi^{\alpha-1}_{t,u}(x)\over dx}\right| \, dxdtdu.
\end{eqnarray*}

We now need to split the proof depending on $\alpha$. First, in the case that $0< \alpha \le 1$, $(\varphi^{\alpha-1}_{t,u})'(x)$ is nonnegative  and
$$
 \int_0^t x^\alpha {d\varphi^{\alpha-1}_{t,u}(x)\over dx}dx\le t\int_0^t x^{\alpha-1} {d\varphi^{\alpha-1}_{t,u}(x)\over dx}dx={t^\alpha\over \Gamma(\alpha)}\left((u-t)^{\alpha-1}-u^{\alpha-1}\right)\le {t^\alpha (u-t)^{\alpha-1}\over \Gamma(\alpha)}
$$
for $ t<u$, where we have applied \cite[Lemma 2.2]{[GP]};  similarly
$$
 \int_0^u x^\alpha {d\varphi^{\alpha-1}_{t,u}(x)\over dx}dx\le {u^\alpha (t-u)^{\alpha-1}\over \Gamma(\alpha)},
$$
for $u<t$. Then we conclude that

\begin{eqnarray*}
A_\alpha(f,g) &=&\int_0^\infty \vert W^{\alpha}_+g(u)\vert  \int_0^u \vert W^{\alpha}_+f(t)\vert  \int_0^t x^\alpha \left|{d\varphi^{\alpha-1}_{t,u}(x)\over dx}\right| \, dxdtdu\\
&\,&\qquad +
\int_0^\infty \vert W^{\alpha}_+g(u)\vert  \int_u^\infty \vert W^{\alpha}_+f(t)\vert  \int_0^u x^\alpha \left|{d\varphi^{\alpha-1}_{t,u}(x)\over dx}\right| \, dxdtdu\\
&\,& \leq {1\over \Gamma(\alpha)}\int_0^\infty \vert W^{\alpha}_+g(u)\vert  \int_0^u \vert W^{\alpha}_+f(t)\vert t^\alpha (u-t)^{\alpha-1}dtdu\\
&\,&\qquad+{1\over \Gamma(\alpha)}\int_0^\infty u^\alpha \vert W^{\alpha}_+g(u)\vert  \int_u^\infty \vert W^{\alpha}_+f(t)\vert  (t-u)^{\alpha-1}dtdu\\
&\,&={1\over \Gamma(\alpha)}\int_0^\infty t^\alpha\vert W^{\alpha}_+f(t)\vert  \int_t^\infty (u-t)^{\alpha-1}\vert W^{\alpha}_+g(u)\vert  dudt\\
&\,&\qquad+{1\over \Gamma(\alpha)}\int_0^\infty u^\alpha \vert W^{\alpha}_+g(u)\vert  \int_u^\infty  (t-u)^{\alpha-1}\vert W^{\alpha}_+f(t)\vert  dtdu\\
&\,&=\int_0^\infty t^\alpha \vert W^{\alpha}_+f(t)\vert \widetilde{g}(t)dt+ \int_0^\infty u^\alpha \vert W^{\alpha}_+g(u)\vert  \widetilde{f}(u)du,
\end{eqnarray*}
where $\widetilde{f}(u):=\displaystyle{{1\over \Gamma(\alpha)}\int_u^\infty (t-u)^{\alpha-1}\vert W^\alpha_+ f (t)\vert dt}$ for $u>0$. Note that $\widetilde{f} \in \Tpa$, with $\|\widetilde{f}\|_{\alpha,p} = \|f\|_{\alpha,p}$. Then, by the same reasoning as in \ref{Holderreason}, we have that
$$
\int_0^\infty t^\alpha \vert W^{\alpha}_+f(t)\vert \widetilde{g}(t)dt\le \Gamma(\alpha+1)\Vert  f\Vert_{\alpha,p}\Vert \widetilde{g}\Vert_{0,p'}\le C_\alpha\Vert  f\Vert_{\alpha,p}\Vert \widetilde{g}\Vert_{\alpha,p'}=C_\alpha\Vert  f\Vert_{\alpha,p}\Vert {g}\Vert_{\alpha,p'}.
$$

So, if $0<\alpha \leq 1$, the proof is finished. In the case that $ \alpha \ge 1$, the function $(\varphi^{\alpha-1}_{t,u})'(x)$ is nonpositive  and
$$
 -\int_0^t x^\alpha {d\varphi^{\alpha-1}_{t,u}(x)\over dx}dx \leq {t^\alpha\over \Gamma(\alpha)}\left(u^{\alpha-1}- (u-t)^{\alpha-1}\right)\le {t^\alpha (u-t)^{\alpha-1}\over \Gamma(\alpha)}\left( \left(1+{t\over u-t}\right)^{\alpha-1}-1\right)
$$
for $ t<u$ (\cite[Lemma 2.2]{[GP]});  similarly
$$
 -\int_0^u x^\alpha {d\varphi^{\alpha-1}_{t,u}(x)\over dx}dx\le {u^\alpha (t-u)^{\alpha-1}\over \Gamma(\alpha)}\left( \left(1+{u\over t-u}\right)^{\alpha-1}-1\right)
$$
for $u<t$. Now we consider the following inequalities,
\begin{eqnarray*}
(1+x)^\nu-1 &\le& \nu x, \qquad \qquad \qquad \: \: \: \: x>0,\,\, 0\le \nu\le 1, \\
(1+x)^\nu-1 &\le& 2^{\nu-1}(\nu  x+x^{\nu}), \qquad  x>0,\,\, 1\le \nu,
\end{eqnarray*}
to conclude that, by the first one, for $1\le \alpha \le 2$
\begin{eqnarray*}
A_\alpha(f,g) &=&\int_0^\infty \vert W^{\alpha}_+g(u)\vert  \int_0^u \vert W^{\alpha}_+f(t)\vert  \int_0^t x^\alpha \left|{d\varphi^{\alpha-1}_{t,u}(x)\over dx}\right| \, dxdtdu\\
&\,&\qquad +
\int_0^\infty \vert W^{\alpha}_+g(u)\vert  \int_u^\infty \vert W^{\alpha}_+f(t)\vert  \int_0^u x^\alpha \left|{d\varphi^{\alpha-1}_{t,u}(x)\over dx}\right| \, dxdtdu\\
&\,&\leq {1\over \Gamma(\alpha-1)}\int_0^\infty t^{\alpha+1}\vert   W^{\alpha}_+f(t)\vert   \int_t^\infty(u-t)^{\alpha-2}\vert W^{\alpha}_+g(u)\vert  dudt\\
&\,&\qquad+{1\over \Gamma(\alpha-1)}\int_0^\infty u^{\alpha+1} \vert W^{\alpha}_+g(u)\vert  \int_u^\infty (t-u)^{\alpha-2}\vert W^{\alpha}_+f(t)\vert  dtdu\\
&\,&=\int_0^\infty t^{\alpha+1} \vert W^{\alpha}_+f(t)\vert W^{1}_+\widetilde{g}(t)dt+ \int_0^\infty u^{\alpha+1} \vert W^{\alpha}_+g(u)\vert  W^1_+\widetilde{f}(u)du,\\
&\,&\le \Gamma(\alpha+1)\left(\Vert  f\Vert_{\alpha,p}\Vert \widetilde{g}\Vert_{1,p'}+\Vert  g\Vert_{\alpha,p'}\Vert \widetilde{f}\Vert_{1,p}\right)\le C_\alpha\Vert  f\Vert_{\alpha,p}\Vert {g}\Vert_{\alpha,p'}.
\end{eqnarray*}
Finally for $\alpha \ge 2$, second inequality leads to
\begin{eqnarray*}
A_\alpha(f,g) &=&\int_0^\infty \vert W^{\alpha}_+g(u)\vert  \int_0^u \vert W^{\alpha}_+f(t)\vert  \int_0^t x^\alpha \left|{d\varphi^{\alpha-1}_{t,u}(x)\over dx}\right| \, dxdtdu\\
&\,&\qquad +
\int_0^\infty \vert W^{\alpha}_+g(u)\vert  \int_u^\infty \vert W^{\alpha}_+f(t)\vert  \int_0^u x^\alpha \left|{d\varphi^{\alpha-1}_{t,u}(x)\over dx}\right| \, dxdtdu\\
&\,&\leq {2^{\alpha-1}\over \Gamma(\alpha-1)}\int_0^\infty t^{\alpha+1}\vert   W^{\alpha}_+f(t)\vert   \int_t^\infty(u-t)^{\alpha-2}\vert W^{\alpha}_+g(u)\vert  dudt\\
&\,&\qquad+{2^{\alpha-1}\over \Gamma(\alpha)}\int_0^\infty t^{2\alpha-1}\vert   W^{\alpha}_+f(t)\vert   \int_t^\infty\vert W^{\alpha}_+g(u)\vert  dudt\\
&\,&\qquad+{2^{\alpha-1}\over \Gamma(\alpha-1)}\int_0^\infty u^{\alpha+1} \vert W^{\alpha}_+g(u)\vert  \int_u^\infty (t-u)^{\alpha-2}\vert W^{\alpha}_+f(t)\vert  dtdu\\
&\,&\qquad+{2^{\alpha-1}\over \Gamma(\alpha)}\int_0^\infty u^{2\alpha-1} \vert W^{\alpha}_+g(u)\vert  \int_u^\infty \vert W^{\alpha}_+f(t)\vert  dtdu\\
&\,&\le C_\alpha\left(\int_0^\infty t^{\alpha+1} \vert W^{\alpha}_+f(t)\vert W^{1}_+\widetilde{g}(t)dt+ \int_0^\infty t^{2\alpha-1} \vert W^{\alpha}_+f(t)\vert W^{\alpha-1}_+\widetilde{g}(t)dt\right)\\
&\,&\qquad+ C_\alpha\left( \int_0^\infty u^{\alpha+1} \vert W^{\alpha}_+g(u)\vert  W^1_+\widetilde{f}(u)du+ \int_0^\infty u^{2\alpha-1} \vert W^{\alpha}_+g(u)\vert  W^{\alpha-1}_+\widetilde{f}(u)du\right)\\
&\,&\le C_\alpha\left(\Vert  f\Vert_{\alpha,p}(\Vert \widetilde{g}\Vert_{1,p'}+\Vert \widetilde{g}\Vert_{\alpha-1,p'})+\Vert  g\Vert_{\alpha,p'}(\Vert \widetilde{f}\Vert_{1,p}+\Vert \widetilde{f}\Vert_{\alpha-1,p})\right)\le C_\alpha\Vert  f\Vert_{\alpha,p}\Vert {g}\Vert_{\alpha,p'},
\end{eqnarray*}
and the proof is finished.\edproof

To finish this section, we recall that the family of operators  $(T_{t,p})_{t\in \RR}$ on $\Tpa$, indicated in \ref{isometry}, which are defined by
\begin{equation}\label{grupo}
T_{t,p} f(s):=e^{-\frac{t}{p}}f(e^{-t}s), \qquad f \in
\mathcal{T}_p^{(\alpha)}(t^\alpha), (p\ge 1 \hbox{ and }\alpha\ge 0),
\end{equation}
is a $C_0$-group of isometries on
$\mathcal{T}_p^{(\alpha)}(t^\alpha)$ (\cite[Theorem 2.5]{[LMPS]}). The infinitesimal
generator $\Lambda$ is given by
$$
(\Lambda f)(s):=-sf'(s)-\frac{1}{p}f(s)
$$
with domain $D(\Lambda)=
\mathcal{T}_p^{(\alpha+1)}(t^{\alpha+1})$;  the point spectrum $\sigma_p(\Lambda)=\emptyset;$ and the usual $\sigma(\Lambda)=i\mathbb{R}$ (\cite[Proposition 2.6]{[LMPS]}). Finally
the semigroups $\{T_{t,p}\}_{t\geq 0}$ and $\{T_{-t,p'}\}_{t\geq
0}$ are adjoint operators of each other acting on
$\mathcal{T}_p^{(\alpha)}(t^\alpha)$ and
$\mathcal{T}_{p'}^{(\alpha)}(t^\alpha)$ with  ${1\over p}+ {1\over
p'}=1$ (\cite[Proposition 2.7]{[LMPS]}).

As it is commented in the Introduction, the $C_0$-group $(T_{t,p})_{t\in \RR}$ will be key to subordinate the generalized Stieltjes operator $\Sbmu$ on the spaces $\mathcal{T}_p^{(\alpha)}(t^{\alpha})$, which will be done in the next section.

\section{Generalized Stieltjes operators on Sobolev spaces}

\setcounter{theorem}{0} \setcounter{equation}{0}

For $\beta, \mu\in \R$, the generalized Stieltjes operator $\mathcal{S}_{\beta, \mu} $ on $\R^+$ is  defined  by
$$
\mathcal{S}_{\beta, \mu} f(t):={t^{\mu-\beta}}\int_0^\infty
{s^{\beta-1}\over (s+t)^{\mu}}f(s)ds=\int_0^\infty
{u^{\beta-1}\over (1+u)^{\mu}}f(tu)du , \quad
\,\,t> 0
$$
for functions defined on $\R^+$, whenever the expression above may be well defined. First of all, we check how these operators act on some particular functions.

\begin{example}
{\rm (i) For $\gamma>0$, we define  $g_\gamma$  by $$
g_{\gamma}(t):= \frac{t^{\gamma-1}}{\Gamma(\gamma)}, \quad t>0.
$$
Then $\mathcal{S}_{\beta, \mu} (g_\gamma)=B(\beta+\gamma-1,\mu-\beta-\gamma+1)g_\gamma $ for $\mu>\beta +\gamma-1>0$. Under these conditions, functions $g_\gamma$ are eigenfunctions of $\mathcal{S}_{\beta, \mu}$. In particular $\mathcal{S}_{\beta, \mu} (\chi_{(0, \infty)})=B(\beta,\mu-\beta)\chi_{(0, \infty)} $ for $0<\beta<\mu$. \\
\noindent (ii) For $\rho + \mu > \beta > 0$, we get
\begin{equation}\label{21}
\mathcal{S}_{\beta, \mu} \left({1\over (1+s)^{\rho}}\right)(t)= B(\beta,\rho+\mu-\beta)\,_2F_1(\rho,\beta;\rho+\mu;1-t), \qquad t>0,
\end{equation}
\cite[p. 317, 3197(1)]{[GR]}; in particular, for $t>0$, we have that
\begin{eqnarray*}
\mathcal{S}_{1, \rho+1} ( (1+s)^{\rho-1})(t)&=&{t^\rho-1\over \rho(t-1)}, \qquad \rho\not=0.\\
\mathcal{S}_{1, 1} \left({1\over 1+s}\right)(t)&=&{\log(t)\over t-1},\\
\mathcal{S}_{\beta, \beta} \left({1\over (1+s)^{\beta}}\right)(t)&=&B(\beta,\beta)\,_2F_1(\beta,\beta;2\beta;1-t),\\
\mathcal{S}_{{1\over 2}, {1\over 2}} \left({1\over \sqrt{1+s}}\right)(t)&=&
\begin{cases}
	2 K (\sqrt{1-t}), & \text{if} \:\: 0<t<1
	\\
	\frac{2}{\sqrt{t}}K\left(\sqrt{1-\frac{1}{t}}\right), & \text{if} \:\: t>1.
\end{cases}
\end{eqnarray*}
 where the function $K(\cdot)$ stands for the complete elliptic integral of the first kind, see for example \cite[p. 255, 3.131(8)]{[GR]}.

\noindent (iii) Take $e_\lambda(s):=e^{-\lambda s}$ for $s>0$ and
$\lambda \in \CC^+$. For $\beta, \mu >0$, we have that
$$
\mathcal{S}_{\beta, \mu} (e_\lambda)(t)= {\mathcal L}\left({s^{\beta-1}\over(1+s)^{\mu}}\right)(\lambda t)
\qquad t>0.
$$
In particular for $\mu\not \in \NN$ and $\beta>0$,
\begin{eqnarray*}
\mathcal{S}_{\beta, \mu} (e_\lambda)(t)&=&{\pi^2\over (\lambda t)^\beta \Gamma(\mu)\sin(\pi(\beta-\mu))}\left({ (\lambda t)^\mu L_{-\mu}^{\mu-\beta}(\lambda t)\over \sin(\pi\mu)\Gamma(1-\beta)}-{ (\lambda t)^\beta L_{-\beta}^{\beta-\mu}(\lambda t)\over \sin(\pi\beta)\Gamma(1-\mu)}\right); \\
\mathcal{S}_{{1\over 2}, {1\over 2}} (e_\lambda)(t)&=&{e^{\frac{\lambda t}{2}}}K_0\left(\frac{\lambda t}{2}\right); \\
\mathcal{S}_{\beta, 1} (e_\lambda)(t)&=&e^{\lambda t}\Gamma(\beta)\Gamma(-\beta+1,\lambda t);\\
\mathcal{S}_{\beta, \beta+{1\over 2}} (e_\lambda)(t)&=&2^\beta \Gamma(\beta) e^{\lambda t\over 2}D_{-2\beta}(\sqrt{2\lambda t});\\
\mathcal{S}_{\beta, \beta-{1\over 2}} (e_\lambda)(t)&=&2^{\beta-{1\over 2}} \Gamma(\beta) {e^{\lambda t\over 2}\over \sqrt{\lambda t}}D_{1-2\beta}(\sqrt{2\lambda t});
\end{eqnarray*}
where these equalities and special functions $L_{-\mu}^{\mu-\beta},\, K_0, \, D_{-2\beta}$ and  the incomplete Gamma function $\Gamma(\cdot, \cdot)$  may be found in \cite[p. 348, 3.383]{[GR]}.
}
\end{example}

As we have commented in the Introduction, the operator $\mathcal{S}_{\beta,\mu}$ defines a bounded operator on $L^p(\R^+)$ for $0<\beta-{1\over p}<\mu$ and $p\ge 1$, obtaining that
$
\Vert \mathcal{S}_{\beta, \mu}\Vert= B( \beta-1/p,\mu-\beta+1/p)
$
; in particular, $\Vert\mathcal{S}\Vert={\pi\over \sin(\pi p)}$ (\cite[Section 9.5, p. 232]{Ha-Li-Po-64}). To extend  this result
in the family of spaces $\mathcal{T}_p^{(\alpha)}(t^\alpha)$  for $p\ge 1$ and $\alpha\ge 0$, we need the next lemma which shows a key commutativity property.

\begin{lemma}\label{llave}
Take $\alpha \ge 0 $ and $\beta, \mu \in \R$. Then $D^\alpha_+
\mathcal{S}_{\beta,\mu}=\mathcal{S}_{\beta,\mu} D^\alpha_+$, i.e.,
$$
D^\alpha_+(\mathcal{S}_{\beta,\mu}(f))=\mathcal{S}_{\beta,\mu}(
D^\alpha_+(f)),\qquad f\in {\mathcal S}_+,
$$
where $D^\alpha_+(t)=\displaystyle{1\over
\Gamma(\alpha+1)}t^\alpha W^\alpha_+ f(t)$ for $t>0$ and $f\in {\mathcal
S}_+$.
\end{lemma}

\bgproof By the equality (\ref{escala}), we have that
\begin{eqnarray*}
\mathcal{S}_{\beta,\mu}(
D^\alpha_+(f))(t)&=& \displaystyle{1\over
\Gamma(\alpha+1)}\int_0^\infty
{u^{\beta-1}\over (1+u)^{\mu}}(tu)^\alpha
W^\alpha_+f(tu)du \cr &=&\displaystyle{1\over
\Gamma(\alpha+1)}t^\alpha W^\alpha_+\left(\int_0^\infty
{u^{\beta-1}\over (1+u)^{\mu}}f(tu)du
\right)(t)=D^\alpha_+(\mathcal{S}_{\beta,\mu}(f))(t),
\end{eqnarray*}
for $t>0$,  $f\in {\mathcal S}_+$ and we conclude the proof. \edproof

The first main result in this section is the following theorem.

\begin{theorem}\label{cotaTa}
The operator $\mathcal{S}_{\beta, \mu}$ is a bounded operator on
$\mathcal{T}_p^{(\alpha)}(t^\alpha)$ and
$$
||\mathcal{S}_{\beta, \mu}||=B(\mu - \beta +1/p, \beta -1/p),
$$
for $ 0< \beta-1/p < \mu$ and $p\ge 1$. If $f\in
\mathcal{T}_p^{(\alpha)}(t^\alpha)$, then
\begin{equation}\label{integral}
\mathcal{S}_{\beta,\mu} f(t)=\displaystyle\int_{-\infty}^\infty
\phi_{\mu -\beta +1/p, \mu}(r)T_{r,p}f(t)dr, \quad t\geq 0,
\end{equation}
where the group $(T_{r, p})_{r \in \R}$ is defined in (\ref{grupo}), and the set of functions $\phi_{\beta,\mu}$ in (\ref{functionp}).
\end{theorem}

\bgproof Let $0< \beta - 1/p < \mu$, $1 \le p$ and $f\in \mathcal{T}_p^{(\alpha)}(t^\alpha)$ be given.  We apply the change of variable $s= t e^{-r}$ to get that
$$
\mathcal{S}_{\beta,\mu} f(t) := t^{\mu - \beta}\int_0^\infty
{\frac{s^{\beta-1}}{(t+s)^\mu}}f(s)ds = \int_{-\infty}^\infty
{\frac{e^{(\mu -\beta+1/p)r}}{(1+e^r)^\mu}}e^{-r/p}f(te^{-r})dr,
$$
and the equality (\ref{integral}) is proved. Observe that by this
equality, the operator $\mathcal{S}_{\beta,\mu}$ is well defined  and is a bounded
operator on $\mathcal{T}_p^{(\alpha)}(t^\alpha)$ for $p\ge 1$. Indeed, we have
\begin{eqnarray*}
||\mathcal{S}_{\beta,\mu} f||_{\alpha, p} &\leq& \int_{-\infty}^\infty ||\phi_{\mu - \beta +1/p,\mu}(r)T_{r,p}f||_{\alpha,p}dr \\
&=&\int_{-\infty}^\infty \phi_{\mu - \beta +1/p,\mu}(r)||f||_{\alpha,p}dr =
B(\mu - \beta +1/p, \beta -1/p) ||f||_{\alpha,p},
\end{eqnarray*}
where we apply Theorem \ref{cottas} (ii). To check the exact value of $||\mathcal{S}_{\beta, \mu} ||$, we apply Lemma \ref{llave},  the
isometry $D_+^{\alpha}$ (see Section 3) and the boundedness of $\mathcal{S}_{\beta,\mu}$ on
$L^p(\RR^+)$, to get
$$
\begin{array}{lcl}
\|\mathcal{S}_{\beta, \mu}\| &=& \displaystyle\sup_{f\neq 0}
\frac{\| \mathcal{S}_{\beta, \mu} f\|_{\alpha,p}}{\|f\|_{\alpha,p}} = \displaystyle \sup_{f\neq 0} \frac{\| D^{\alpha}_+
\mathcal{S}_{\beta, \mu} f\|_{p}}{\|D^{\alpha}_+ f\|_{p}} \\
\\ &=& \displaystyle \sup_{f\neq 0} \frac{\|
\mathcal{S}_{\beta, \mu}  D^{\alpha}_+  f\|_{p}}{\|D^{\alpha}_+
f\|_{p}} = \displaystyle\sup_{g\neq 0} \frac{\|
\mathcal{S}_{\beta, \mu} g\|_{p}}{\|g\|_{p}} = B(\mu - \beta +1/p, \beta -1/p).
\end{array}
$$
Finally we conclude the proof.
\edproof

\begin{remark} {\rm  In the case $p=1$ and $\beta=\mu=1$ we remark that the Stieltjes operator $\mathcal{S}$ does not
take $\mathcal{T}_1^{(\alpha)}(t^\alpha)$ in
$\mathcal{T}_1^{(\alpha)}(t^\alpha)$. Hence,  the function $h_2$,  given by
$h_2(t):=(1+t)^{-2}$ for $t>0$, belongs to
$\mathcal{T}_1^{(\alpha)}(t^\alpha)$ and
$$
\mathcal{S}
h_2(t)={t-\log(t)-1\over (t-1)^2}, \qquad t>0.
$$
Since $\mathcal{S}
h_2$ does not belong to
$L^1(\mathbb{R}^+)$ and
$\mathcal{T}_1^{(\alpha)}(t^{\alpha})\hookrightarrow
L^1(\mathbb{R}^+)$, we obtain that
$\mathcal{S}
h_2\not \in
\mathcal{T}_1^{(\alpha)}(t^\alpha)$.

}
\end{remark}

In the next result, we are able to describe
$\sigma(\mathcal{S}_{\beta, \mu})$ for suitable $\beta, \mu>0$.

\begin{theorem}\label{spec}
Let $1\le p<\infty,$ and $\mathcal{S}_{\beta,\mu}
:\mathcal{T}_p^{(\alpha)}(t^\alpha) \to
\mathcal{T}_p^{(\alpha)}(t^\alpha)$ the generalized Stieltjes
operator with $0< \beta -1/p< \mu$. Then

$$
\sigma(\mathcal{S}_{\beta,\mu})= {\biggl\{ B\left(\beta-{1\over p}+i\xi,\mu-\beta+{1\over p}-i\xi\right)\,:\, \xi \in \R \biggr\}}\cup\{0\}.
$$
\end{theorem}

\bgproof  Recall $(T_{t,p})_{t\in \RR}$ is an uniformly bounded
$C_0$-group  whose infinitesimal
generator is $(\Lambda, D(\Lambda))$ and $\mathcal{S}_{\beta,\mu}=\theta(\phi_{\mu - \beta +1/p,\mu})$, see  Theorem
\ref{cotaTa}. By  \cite[Theorem 3.1]{[Se]}, we obtain
$$
\sigma(\mathcal{S}_{\beta,\mu})=\overline{\widehat{\phi_{\mu - \beta +1/p,\mu}}(\sigma(i\Lambda))}.
$$
 As $\sigma(i\Lambda)=\RR$ (see
Section 3), we apply Theorem \ref{cottas} (iii) to conclude
$$
\sigma(\mathcal{S}_{\beta,\mu})= {\biggl\{ B\left(\beta-{1\over p}+i\xi,\mu-\beta+{1\over p}-i\xi\right)\,:\, \xi\in\R \biggr\}}\cup\{0\},
$$
where we have applied that $\lim_{\xi\to\pm\infty}\Gamma(a+i\xi)=0$ for $a>0$.\edproof

\begin{remark}
{\rm  In the case that $\mu=1$ and  $0<\beta- {1\over p}<1$,  we obtain that
\begin{eqnarray*}
\sigma(\mathcal{S}_{\beta, 1})%&=& \left\{w\in\mathbb{C}: \left|\left(\frac{w}{\Gamma(\beta+1)}\right)^{1/\beta}-\frac{p}{2(p-1)}\right|=\frac{p}{2(p-1)}\right\}\\
&=&{\left\{{\pi\over \sin(\pi(\beta-{1\over p}+i\xi))}: \quad \xi \in \R\right\}}\cup\{0\},
\end{eqnarray*}
where we apply the   Euler's reflection formula $\Gamma(z)\Gamma(1-z)=\displaystyle{\pi\over \sin(\pi z)}$,  $z\not \in \ZZ$. For $\beta=1$ and $p=2$, we have that
$$
\sigma(\mathcal{S}_{1,1})= {\left\{{\pi\over \cosh(\pi\xi)}: \quad \xi \in \R\right\}}\cup\{0\}=[0,\pi].
$$
In the last section we draw some of these families  of spectrum sets.}
\end{remark}

Now we identify the adjoint generalized Stieltjes operator
$(\mathcal{S}_{\beta,\mu})'$  on $\mathcal{T}_{p'}^{(\alpha)}(t^\alpha)$.

\begin{theorem}\label{dualStieltjes} For $p>1$ and  $0< \beta -1/p < \mu$, the adjoint generalized Stieltjes operator of
$\mathcal{S}_{\beta,\mu}$  on
$\mathcal{T}_p^{(\alpha)}(t^\alpha)$ is $\mathcal{S}_{\mu-\beta+1,\mu} $  on $\mathcal{T}_{p'}^{(\alpha)}(t^\alpha)$,
i.e.
$$
\langle \mathcal{S}_{\beta,\mu} f,g \rangle_\alpha=\langle
f,\mathcal{S}_{\mu-\beta+1,\mu}g \rangle_\alpha, \qquad f \in
\mathcal{T}_p^{(\alpha)}(t^\alpha), \quad  g
\in\mathcal{T}_{p'}^{(\alpha)}(t^\alpha),
$$
where $\langle \quad,\quad \rangle_\alpha$ is given in Section 3 and $\frac{1}{p}+\frac{1}{p'}=1$. Note that
$\mathcal{S}_{\beta, \mu}$ is an injective, non-surjective and of dense range  on $\mathcal{T}_p^{(\alpha)}(t^\alpha)$.
\end{theorem}

\bgproof  First, as the operator $\mathcal{S}_{\beta,\mu}$ intertwines with the operator $D^\alpha_+$ (Lemma \ref{llave}), we apply the Fubini theorem to get that
\begin{eqnarray*}
\langle \mathcal{S}_{\beta,\mu}f, g\rangle_\alpha &=&\langle D^\alpha_+ \mathcal{S}_{\beta,\mu},D^\alpha_+ g
\rangle_{0}= \langle  \mathcal{S}_{\beta,\mu} D^\alpha_+ f, D^\alpha_+ g\rangle_0 \\
&=&
\int_0^\infty   \int_0^\infty \frac{t^{\mu -\beta} s^{\beta -1}}{(t+s)^\mu}   D^\alpha_+ g(t) D^\alpha_+ f (s) ds dt
\\
&=& \langle   D^\alpha_+ f,  \mathcal{S}_{\mu -\beta +1,\mu} D^\alpha_+ g\rangle_0
= \langle f,\mathcal{S}_{\mu -\beta +1,\mu} g\rangle_\alpha.
\end{eqnarray*}
Secondly, the injectivity of the Laplace transform $\mathcal{L}$  and (\ref{faclap}) implies  the injectivity of $\mathcal{S}_{\beta,\mu}$ on $L^p(\R^+)$ and then in $\mathcal{T}_p^{(\alpha)}(t^\alpha) \subseteq L^p(\R^+)$. By Theorem \ref{spec}, $0\in \sigma(\mathcal{S}_{\beta,\mu})$ and $\mathcal{S}_{\beta,\mu}$ is not invertible so, by the open mapping theorem, $\mathcal{S}_{\beta,\mu}$ cannot be surjective. Moreover, as the adjoint of $\mathcal{S}_{\beta,\mu}$ is injective since it is another generalized Stieltjes operator, we conclude that $\mathcal{S}_{\beta,\mu}$ is of dense range on $\mathcal{T}_p^{(\alpha)}(t^\alpha)$.
\edproof

Then, the fact that generalized Stieltjes operators commute, immediately leads us to the following  corollary.

\begin{corollary}\label{dualll}
	Let $\mathcal{S}_{\beta,\mu}$ and $\mathcal{S}_{\gamma,\nu}$ be generalized Stieltjes operators on $\mathcal{T}_{p}^{(\alpha)}(t^\alpha)$ and $\mathcal{T}_{p'}^{(\alpha)}(t^\alpha)$ respectively for $p>1$. Then, $\mathcal{S}_{\beta,\mu} (\mathcal{S}_{\gamma,\nu})' = (\mathcal{S}_{\gamma,\nu})'  \mathcal{S}_{\beta,\mu}$. As a consequence, $\mathcal{S}_{\beta,\mu}$ is a normal operator on $\mathcal{T}_{2}^{(\alpha)}(t^\alpha)$, while $\mathcal{S}_{\beta,2\beta-1}$ is also self-adjoint for $\beta>{1\over 2}.$
\end{corollary}

\begin{remark}\label{autoad} {\rm As  the operator $\mathcal{S}_{\beta,2\beta-1}$ is  self-adjoint on $\mathcal{T}_{2}^{(\alpha)}(t^\alpha)$ for $\beta>{1\over 2},$ then the spectrum $\sigma(\mathcal{S}_{\beta,2\beta-1})$ is a subset of real numers,
\begin{eqnarray*}
	\sigma(\mathcal{S}_{\beta,2\beta-1})&=& {1\over \Gamma(2\beta-1)}{\left\{\Gamma( \beta -1/2+i\xi)\Gamma(\beta -1/2-i\xi): \quad \xi \in \R\right\}}\cup\{0\}\\&=&[0,B(\beta -1/2, \beta -1/2)],
\end{eqnarray*}
where we have used that $\Gamma(z)\Gamma(\overline{z})\in \R$.  This result was proved in $L^2(\R^+)$ for the operator ${\mathcal S}$ in \cite[p. 169]{Ca}  and  finally for $\mathcal{S}_{\beta,2\beta-1}$ for $\beta>{1\over 2}$ in  \cite[Proposition 1.1]{Federer}.}
\end{remark}

To finish this section, we give an explicit formula for the composition of two generalized Stieltjes operators involving the Gaussian hypergeometric function $_2F_1$.

\begin{theorem} \label{compConDual2} Let $\mathcal{S}_{\beta, \mu}$ and $\mathcal{S}_{\gamma, \nu}$ be the generalized Stieltjes operators on $\mathcal{T}_p^{(\alpha)}(t^\alpha)$ for $p\ge 1, \: \alpha \geq 0$. Then
\begin{eqnarray*}
\frac{\mathcal{S}_{\beta, \mu}  \mathcal{S}_{\gamma, \nu}f(t)}{B(\beta +\nu -\gamma, \gamma +\mu -\beta)} =
\int_0^t f(s) \frac{s^{\beta -1}}{t^\beta}{}_2 F_1\left(\mu, \beta +\nu -\gamma; \mu +\nu; 1-{s \over t}\right)ds
\\
+\int_t^\infty f(s)\frac{t^{\mu -\beta}}{s^{\mu - \beta +1}} {}_2 F_1\left(\mu, \gamma +\mu -\beta ; \mu +\nu;
1-{t \over s}\right)ds
\end{eqnarray*}
In particular, if $p>1$ and $\mathcal{S}_{\gamma, \nu}$ is a generalized Stieltjes operator on the dual space $\mathcal{T}_{p'}^{(\alpha)}(t^\alpha)$, where ${1 \over p} + {1 \over {p'}} = 1$, we have that
\begin{eqnarray*}
\frac{\mathcal{S}_{\beta, \mu}  (\mathcal{S}_{\gamma, \nu})' f(t)}{B(\beta +\gamma -1, \nu -\gamma +\mu -\beta +1)} =
\int_0^t f(s) \frac{s^{\beta -1}}{t^\beta}{}_2 F_1\left(\mu, \beta +\gamma -1; \mu +\nu; 1-{s \over t}\right)ds
\\
+\int_t^\infty f(s)\frac{t^{\mu -\beta}}{s^{\mu - \beta +1}} {}_2 F_1\left(\mu, \nu -\gamma +\mu -\beta +1; \mu +\nu; 1-{t \over s}\right)ds
\end{eqnarray*}
for $f\in \mathcal{T}_p^{(\alpha)}(t^\alpha)$ and $t$ almost everywhere on $\R^{+}$.
\end{theorem}

\bgproof Describing the generalized Stieltjes operators with the Sinclair map (see (\ref{integral})) and the functions $\phi_{\beta, \mu}$ given in (\ref{functionp}), it follows that $\mathcal{S}_{\beta, \mu}  \mathcal{S}_{\gamma, \nu} = \theta(\phi_{\mu -\beta +1/p, \mu} * \phi_{\nu - \gamma + 1/p, \nu})$. This convolution has already been obtained in Theorem \ref{cottas} (iv). Then, applying the change of variable $te^{-r} = s$ one concludes that
\begin{eqnarray*}
\frac{\mathcal{S}_{\beta, \mu}  \mathcal{S}_{\gamma, \nu}f(t)}{B(\beta +\nu -\gamma, \gamma +\mu -\beta)} &=&
\int_{-\infty}^\infty e^{-r(\beta -1/p)} {}_2 F_1(\mu, \beta +\nu -\gamma ;\mu +\nu;  1-e^{-r}) T_{r,p}(f)(t) dr
\\ &=& \int_{0}^\infty f(s) \frac{s^{\beta -1}}{t^\beta} {}_2 F_1\left(\mu, \beta +\nu -\gamma ;\mu +\nu;  1-{s \over t}\right)ds.
\end{eqnarray*}
Then, one can split the domain $(0, \infty)$ into $(0,t)$ and $(t, \infty)$, and then apply the formula (\ref{kummer}) to ${}_2 F_1(\cdot)$ in the integral over $(t, \infty)$, obtaining the first formula of the theorem.
\\
The second formula, involving $(\mathcal{S}_{\gamma, \nu})'$, is obtained directly from the first one and applying the expression for its adjoint operator given in Theorem \ref{dualStieltjes}.
\edproof

\section{Composition of generalized Stieltjes, Ces\`{a}ro operators and Hilbert transform}

\setcounter{theorem}{0} \setcounter{equation}{0}

In this section, we recall the definitions of the Ces\'aro operators and the Hilbert transform, and give a handful of results involving both of them and the generalized Stieltjes operators on the family of $\Tpa$ spaces. Let's begin with the first one.

In \cite{[LMPS]} other family of bounded operators of $\mathcal{T}_p^{(\alpha)}(t^\alpha)$, namely the generalized Ces\`{a}ro operators $\mathcal{C}_\gamma$, given by
$$
\mathcal{C}_\gamma f(t) := \frac{\gamma}{t^\gamma}\int_0^t (t-s)^{\gamma-1}f(s)ds,  \qquad t>0,
$$
with  $\gamma > 0$ are considered  in detail. They  define  bounded operators on $\mathcal{T}_p^{(\alpha)}(t^\alpha)$,
$$
||\mathcal{C}_\gamma||=\frac{\Gamma(\gamma+1)\Gamma(1-1/p)}{\Gamma(\gamma+1-1/p)}=\gamma B(\gamma, 1-1/p),
$$
for $p>1$ and $\alpha \geq 0$. Similarly to the Stieltjes operators, a key tool on the proof is to subordinate Ces\`{a}ro operators in terms of the $C_0$-group $(T_{t,p})_{t\in \RR}$ and write
$$
\mathcal{C}_\gamma f(t)= \gamma\int^{\infty}_0
(1-e^{-r})^{\gamma-1}e^{-r(1-1/p)}T_{r,p}f(t)dr=\theta(\psi_{\gamma,1-1/p})f(t), \quad t\geq
0,
$$
(\cite[Theorem 3.3]{[LMPS]}), where $\psi_{\gamma, \nu}$ are the family of functions given by (\ref{psifunctions}) in Section 2. In \cite[Theorem 3.5]{[LMPS]},  authors also describe the spectrum set of $\mathcal{C}_\gamma$, i.e.,
$$
\sigma(\mathcal{C}_\gamma)=\Gamma(\gamma+1)\overline{\left\{{\Gamma(1-{1\over
p}+i\xi)\over \Gamma(\gamma+1-{1\over p}+i\xi)} \ : \ \xi\in \R
\right\}}= {\biggl\{ \gamma B\left(\gamma,1-{1\over p}+i\xi\right)\,:\, \xi\in\R \biggr\}}\cup\{0\}.
$$

Then, after this short introduction to  generalized Ces\`{a}ro operators, we give a handful of results concerning the composition of these generalized Ces\`{a}ro operators with the generalized Stieltjes operators. In some cases, note that the generalized Ces\`{a}ro operators factorize the Stieltjes operators.

\begin{theorem}\label{CesStieltComp}
Let $0 <\beta -1/p < \mu$, $\gamma >0$ and $1<p$. Then $\mathcal{S}_{\beta, \mu} \mathcal{C}_\gamma = \mathcal{C}_\gamma  \mathcal{S}_{\beta, \mu}$ and

\begin{eqnarray*}
\mathcal{S}_{\beta, \mu}  \mathcal{C}_\gamma f(t) = \gamma B(\gamma, \mu -\beta +1)
t^{\mu -\beta} \int_0^{\infty}\frac{s^{\beta -1}}{(s+t)^{\mu}} {}_2F_1\left(\mu, \gamma;\mu -\beta +1 +\gamma; \frac{t}{s+t}\right) f(s)ds
\end{eqnarray*}
for $t>0$ and $f\in\mathcal{T}_p^{(\alpha)}(t^\alpha)$. In particular,
\begin{itemize}
\item[(i)]for $\mu > \gamma +1-1/p$, the following factorization holds
\begin{eqnarray*}
\mathcal{S}_{\gamma +1, \mu}  \mathcal{C}_{\gamma} = \gamma B(\gamma, \mu - \gamma) \mathcal{S}_{1, \mu - \gamma}.
\end{eqnarray*}
\item[(ii)] for $\beta>{1\over p}$,
\begin{eqnarray*}
\mathcal{S}_{\beta,\beta} \mathcal{C}_1 (f) (t) &=& \frac{1}{(\beta -1)t}
\int_0^{\infty} \left({(t+s)^{\beta-1}-s^{\beta-1}\over(t+s)^{\beta-1}}  \right) f(s)ds,\qquad \beta\not=1,\\
\mathcal{S}_{1,1}  \mathcal{C}_1 (f) (t) &=& \frac{1}{t}
\int_0^{\infty} \log(1+{t\over s}) f(s)ds,
\end{eqnarray*}
for $t>0$ and $f\in\mathcal{T}_p^{(\alpha)}(t^\alpha).$

\end{itemize}

\end{theorem}
\bgproof As $\mathcal{C}_\gamma = \theta(\psi_{\gamma, 1-1/p})$ and $\mathcal{S}_{\beta, \mu} = \theta(\phi_{\mu -\beta +1/p,\mu})$, where $\theta$ is given by (\ref{subor}), we have that
$$
\mathcal{S}_{\beta, \mu} \mathcal{C}_\gamma  =
\theta(\phi_{\mu -\beta +1/p, \mu} * \psi_{\gamma, 1-1/p})=\theta(\psi_{\gamma, 1-1/p}*\phi_{\mu -\beta +1/p, \mu} )=\mathcal{C}_\gamma  \mathcal{S}_{\beta, \mu},
$$
 and both operators commute.

We consider the expression of  $\phi_{\mu -\beta +1/p, \mu} * \psi_{\gamma, 1-1/p}$ given in Proposition \ref{convolu}, to get
\begin{eqnarray*}
&\quad&\mathcal{S}_{\beta, \mu} \mathcal{C}_\gamma f(t) = \theta(\phi_{\mu -\beta +1/p, \mu} * \psi_{\gamma, 1-1/p}) (f) (t)\\
&=& \int_{-\infty}^{\infty}\gamma B(\gamma, \mu -\beta +1) \frac{e^{(\mu -\beta +1/p)r}}{(1+e^r)^\mu}
{}_2F_1\left(\mu, \gamma; \mu -\beta + 1 + \gamma; \frac{e^r}{1+e^r}\right) T_{r,p}f(t)dr
\\ &=& \gamma B(\gamma, \mu -\beta +1)
t^{\mu -\beta} \int_0^{\infty}\frac{s^{\beta -1}}{(s+t)^{\mu}} {}_2F_1\left(\mu, \gamma;\mu -\beta +1 +\gamma; \frac{t}{s+t}\right) f(s)ds,
\end{eqnarray*}
where we have applied  the change of variable $-te^{-r}dr=ds$.

To show the part (i), we apply the identity (\ref{factor}) for $\mu > \gamma +1-1/p$ to get

$$\mathcal{S}_{\gamma +1, \mu} \mathcal{C}_{\gamma} = \theta(\gamma B(\gamma, \mu -\gamma)\phi_{\mu -\gamma -1 +1/p, \mu - \gamma}) = \gamma B(\gamma, \mu - \gamma) \mathcal{S}_{1, \mu - \gamma}.
$$
Finally we apply  Proposition \ref{convolu} (ii) to conclude the proof of the last two equalities. \edproof

Now we turn to a semifinite Hilbert transform, ${\mathcal H_+}: L^p(\R^+)\to L^p(\R^+)$ given by
$$
{\mathcal H}_+f(t):=\hbox{p.v.}{i\over \pi}\int_0^\infty{f(s)\over t-s}ds, \qquad t>0,
$$
for $1<p<\infty$. It is known that the operator ${\mathcal H_+}$ is linear, bounded and   $\displaystyle{\|\mathcal{H}_+\| = \cot \frac{\pi}{2p^\ast}}$
where $p^\ast := \max\{p,p'\}, \frac{1}{p}+\frac{1}{p'}=1,$ see \cite[Section 6.1]{[HKV]}. Then, we give the following theorem about the behaviour of ${\mathcal H_+}$ on $\Tpa$ and with the generalized Stieltjes operator.

\begin{theorem}\label{semihilbert}
Let $p>1$  and  $\alpha \geq 0$ be. Then
\begin{itemize}
\item[(i)] $(D^{\alpha}_+)^{-1} \circ \mathcal{H}_+=\mathcal{H}_+ \circ(D^{\alpha}_+)^{-1} $.
\item[(ii)] $D_{+}^{\alpha} \circ \mathcal{H}_+ = \mathcal{H}_+ \circ D_{+}^{\alpha}$, the operator  $\mathcal{H}_+$ is a bounded on $\mathcal{T}_p^{(\alpha)}(t^\alpha)$, and $$\|\mathcal{H}_+\|= \cot \frac{\pi}{2p^\ast}.$$
\item[(iii)] $\mathcal{H}_+ \circ \mathcal{S}_{\beta, \mu} = \mathcal{S}_{\beta, \mu} \circ \mathcal{H}_+ $ on $\mathcal{T}_p^{(\alpha)}(t^\alpha)$ for  $0 < \beta - 1/p < \mu$.
\end{itemize}
\end{theorem}
\bgproof (i)
Let $f \in L^p(\R^+)$ and $x>0$ be. Then
\begin{eqnarray*}
\mathcal{H}_+( (D^{\alpha}_+)^{-1} f)(x)
&=&\lim_{\varepsilon \rightarrow 0^+}\frac{i}{\pi}\int_{(0,x-\varepsilon)\cup (x + \varepsilon, \infty)} \frac{\alpha }{x-t} \int_1^\infty\frac{(u-1)^{\alpha-1}}{u^\alpha}f(ut)dudt
\\
&=& \lim_{\varepsilon \rightarrow 0^+}\frac{i}{\pi} \alpha \int_1^\infty\frac{(u-1)^{\alpha-1}}{u^\alpha} \int_{(0,x-\varepsilon)\cup (x + \varepsilon, \infty)} \frac{1}{x-t}f(ut)dtdu
\\
&=& \lim_{\varepsilon \rightarrow 0^+} \alpha \int_1^\infty\frac{(u-1)^{\alpha-1}}{u^\alpha} \frac{i}{\pi}
\int_{(0,u(x-\varepsilon))\cup (u(x + \varepsilon), \infty)}
\frac{1}{xu-v}f(v)dvdu.
\end{eqnarray*}
Now, recall that the maximal operator $\mathcal{H}_+^\ast$ defined as $\mathcal{H}_+^\ast f (x) = \sup_{\varepsilon>0}|\mathcal{H}_{+,\varepsilon}f(x)|$ is $(p,p)$-strong (see for example \cite{[Du2]}), where
$$
\mathcal{H}_{+,\varepsilon}f(x):= \frac{i}{\pi}
\int_{(0,u(x-\varepsilon))\cup (u(x + \varepsilon), \infty)}
\frac{1}{xu-v}f(v)dv, \qquad x>0.
$$

 Fixed $x>0$ and $\varepsilon>0$, the following inequality holds
\begin{eqnarray*}
	\left|\frac{(u-1)^{\alpha-1}}{u^\alpha}\mathcal{H}_{+,u\varepsilon}f(xu)\right| \leq \frac{(u-1)^{\alpha-1}}{u^\alpha} \mathcal{H}_+^\ast f(xu), \qquad u>1.
\end{eqnarray*}
As the last function belongs to $L^1(1,\infty)$,  we apply the dominated convergence theorem to conclude the proof.

(ii) By the part (i), we have that $
D_{+}^{\alpha} \circ \mathcal{H}_+ = D_{+}^{\alpha} \circ \mathcal{H}_+ \circ ( (D^{\alpha}_+)^{-1} \circ D^{\alpha}_+)
= \mathcal{H}_+ \circ D^{\alpha}_+
$
To check the norm $\|\mathcal{H}_+\|$ on $\mathcal{T}_p^{(\alpha)}(t^\alpha)$ we have that
$$
	\|\mathcal{H}_+\|=\sup_{f \neq 0} \frac{\|D_+^\alpha\mathcal{H}_+ f\|_{0,p}} {\|D_+^\alpha f\|_{0,p}}
	= \sup_{g \neq 0} \frac{\|\mathcal{H}_+ g\|_{0,p}} {\|g\|_{0,p}}
	 = \cot \frac{\pi}{2p^\ast}.
$$
The proof of the part (iii) is similar to (i) and is left to the reader.
\edproof

To conclude this section, we detail a study of ${\mathcal H_+}$ on $\Tpa$ from the rich point of view of the UMD property.
\begin{remark}
{\rm The geometric property UMD (Unconditional Martingale Difference property) was first introduced by Maurey and Pisier in the study of vector valued martingale theory (\cite{[MaPi]}) in Banach spaces $X$. Later Burkholder together with some other authors developed a rich theory on the UMD property, and in particular  the characterization of the UMD-property in terms of the vector-valued Hilbert transform: let ${\mathcal H}_{\varepsilon, N}$ the following bounded operator defined on $L^p(\R; X)$ for all $\varepsilon \in (0,1)$, $N>1$,
$$
{\mathcal H}_{\varepsilon, N}f(t):={i\over \pi}\int_{\varepsilon\le \vert s\vert\le N}{f(t-s)\over s}ds, \qquad {\hbox{ for a.a. } t\in \R,} \quad f\in L^p(\R; X).
$$
The family $({\mathcal H}_{\varepsilon, N})_{\varepsilon, N}$ admits a strong limit $H$ as $\varepsilon $ goes to $0^+$ and $N$ goes  to $+\infty$ in all $ L^p(\R; X)$, $p\in (1,\infty)$ if and only if $X$ has the UMD-property (\cite{[Bo],[Bu]}).

The classical examples of UMD spaces include all the finite dimensional Banach spaces and any space of the form $L^p(\Omega, \mu)$ for $1 < p < \infty$. If $X$ is any UMD-space, then $X$ is reflexive and $X^*$ is also a
UMD-space; if $Y\subset X$ is a closed subspace then $Y$ is also a UMD-space. See, for example, \cite[Section II]{[Fra86]}. In particular all spaces   $\mathcal{T}_p^{(\alpha)}(  t^\alpha)$ (and $\mathcal{T}_p^{(\alpha)}( \vert t\vert^\alpha)$ considered in the next section) are UMD-spaces for $1<p<\infty$ and $\alpha \ge 0$.

 In UMD-spaces, singular integrals of uniformily bounded $C_0$-groups, $T:(T(t))_{t\in \R}\subset {\mathcal B}(X)$, converges  strongly to a bounded operators, i.e.,
 $$
 \lim_{\varepsilon\to 0^+, \,N \to\infty}{i\over \pi}\int_{\varepsilon\le \vert s\vert\le N}{T(s)x\over s}ds:=H^{T}x, \qquad \hbox{ for } x\in X,
$$
(\cite[Proposition 5.2]{Mo}), which includes the Hilbert transform where $T(s)f=f(\cdot-s)$, for $s\in \R$, and $f\in L^p(\R; X)$. We show how to subordinate also the operator ${\mathcal H}_{+}$  via the $C_0$-group $(T_{r,p})_{r\in \R}$. As a consequence  ${\mathcal H}_{+}$  is a bounded operator on $\mathcal{T}_p^{(\alpha)}(  t^\alpha)$ for $1<p<\infty$ and $\alpha \ge 0$.

  Take $t>0$  and  $0<\varepsilon<1$  and consider operators
$$
{\mathcal H}_{+,\varepsilon}f(t):={i\over \pi}\int_{(0,te^{-\varepsilon})\cup(te^{\varepsilon},\infty)}{f(s)\over t-s}ds, \qquad  {\hbox{ for a.a. } t>0,} \quad f\in \mathcal{T}_p^{(\alpha)}(  t^\alpha).
$$
Now we change the variable $s=te^{-r}$ to obtain
$$
{\mathcal H}_{+,\varepsilon}f(t)={i\over \pi}\int_{\vert r\vert >\varepsilon}{e^{-r}f(te^{-r})\over 1-e^{-r}}dr={i\over \pi}\int_{\vert r\vert >1}{e^{-r\over p'}T_{r,p}f(t)\over 1-e^{-r}}dr+{i\over \pi}\int_{1>\vert r\vert >\varepsilon}{e^{-r\over p'}T_{r,p}f(t)\over 1-e^{-r}}dr
$$
where the $C_0$-group $(T_{r,p})_{r\in\R}$ is defined in (\ref{grupo}). As $(T_{r,p})_{r\in\R}$ is a $C_0$-group of contractions in $\mathcal{T}_p^{(\alpha)}(  t^\alpha)$, the first integral defines a bounded operator. Now note that the second integral is written as
$$
{i\over \pi}\int_{1>\vert r\vert >\varepsilon}{e^{-r\over p'}T_{r,p}f(t)\over 1-e^{-r}}dr={i\over \pi}\int_{1>\vert r\vert >\varepsilon}{(re^{-r\over p'}-1+e^{-r})T_{r,p}f(t)\over (1-e^{-r})r}dr+ {i\over \pi}\int_{1>\vert r\vert >\varepsilon}{T_{r,p}f(t)\over r}dr.
$$
The first summand defines a bounded operator due to the function $h(r)= \displaystyle{re^{-r\over p'}-1+e^{-r}\over (1-e^{-r})r}$ is continuous on $[-1,1]$. As $\mathcal{T}_p^{(\alpha)}(  t^\alpha)$ is a UMD-space, the second summand admits a strong limit in $\mathcal{T}_p^{(\alpha)}(  t^\alpha)$ when $ \varepsilon $ goes to $0^+$ (\cite[Lemma 4.1]{Mo}). The same proof runs on the space $\mathcal{T}_p^{(\alpha)}(  \vert t\vert^\alpha)$, which will be introduced in the next section.

}

\end{remark}

\section{Generalized  Stieltjes operators on $\R$ and Fourier transform}

\setcounter{theorem}{0} \setcounter{equation}{0}

In \cite[Section 4]{[LMPS]}, a family of spaces  $\mathcal{T}_p^{(\alpha)}(|t|^\alpha)$ (which are contained in
$L^p(\RR)$) are presented in a similar way than spaces $\mathcal{T}_p^{(\alpha)}(t^\alpha)$ (embedded into $L^p(\RR^+)$). Now we mention the main properties of these spaces. After that, we will study the generalized Stieltjes operator $\Sbmu$, which first has to be extended to the whole real line $\R$, on these spaces in a similar fashion we have done on $\Tpa$.

Let ${\mathcal S}$ be the Schwartz class on
$\mathbb{R}$ and we set
\begin{eqnarray*}
W_{-}^{-\alpha}f(x)&:=&{1\over
\Gamma(\alpha)}\int_{-\infty}^{x}(x-t)^{\alpha-1}f(t)dt,\cr
W_{-}^{\alpha}f(x)&:=&{1\over\Gamma(n-\alpha)}{d^n \over
dx^n}\int_{-\infty}^{x}(x-t)^{n-\alpha-1}f(t)dt,
\end{eqnarray*}
 $W_{-}^0f:=f$, for $x\in\R$ and a natural number $n>\alpha$.
Putting $\tilde{f}(x)=f(-x)$, it is readily seen that
$W_{+}^{\alpha}f(x) = W_{-}^{\alpha}\tilde{f}(-x)$ for all
$\alpha\in\R$, $f\in {\mathcal S}$ and $x\in\R$. Equalities
$W_{-}^{\alpha+\beta}=W_{-}^\alpha W_{-}^\beta$ and
$W_-^nf=f^{(n)}$ hold for each natural number $n$ and
$\alpha,\beta\in\R$.

For $f\in {\mathcal S}$, put
$$
W^{\alpha}_0f(t):= \left\{\begin{array}{ll}   W^\alpha_-f(t),              & t<0,\\
                                              W^\alpha_+f(t),& t>0.\\
                          \end{array} \right.
$$
For $\lambda >0$, we have that $W_0^\alpha (f_\lambda)=
\lambda^\alpha (W_0^\alpha f)_\lambda$, where $f_\lambda(t):=
f(\lambda t)$ for $t\in \R$.

For $1\le p<\infty$, the Banach space
$\mathcal{T}_p^{(\alpha)}(|t|^\alpha)$  is defined  as the
completion of the Schwartz class on $\mathbb{R}$ in the norm
\begin{equation*}
|||f|||_{\alpha, p}:={1\over \Gamma(
\alpha+1)}\left(\int_{-\infty}^\infty \left(|W^\alpha_0 f(t)\vert
\,\vert t|^{\alpha }\right)^pdt\right)^{\frac{1}{p}},
\end{equation*}
see \cite[Definition 4.1]{[LMPS]}

Properties similar to those of
$\mathcal{T}_p^{(\alpha)}(t^\alpha)$ hold for
$\mathcal{T}_p^{(\alpha)}(|t|^\alpha)$. For example,  the subspace
$\mathcal{T}_1^{(\alpha)}(|t|^\alpha)$  is a subalgebra of
$L^1(\R)$ for the convolution product
\begin{equation}\label{convos2}
f\ast g(t)=\int_{-\infty}^{\infty} f(t-s)g(s)ds, \qquad t\in
\R,\quad  f, g\in \mathcal{T}_1^{(\alpha)}(|t|^\alpha),
\end{equation}
see \cite[Theorem 1.8]{Ga-Mi-06} and also \cite[Theorem 2]{Mi-07} and for $p=2$, the subspace
$\mathcal{T}_2^{(\alpha)}(|t|^\alpha)$  is a Hilbert space, see similar ideas in \cite[Section 2]{[GMMS]}. For  $1< p<\infty$, the Banach space
$\mathcal{T}_p^{(\alpha)}(|t|^\alpha)$  is a module for the
algebra $\mathcal{T}_1^{(\alpha)}(|t|^\alpha)$ (\cite[Theorem 4.3]{[LMPS]}).

Take  $p\geq 1$ and $\beta>\alpha>0$, then
$\mathcal{T}_p^{(\beta)}(\vert t\vert^\beta)\hookrightarrow \mathcal{T}_p^{(\alpha)}(\vert t\vert^\alpha)\hookrightarrow  L^p(\mathbb{R}) $ and the operator $D^\alpha_0: \mathcal{T}_p^{(\alpha)}(\vert t\vert^\alpha)\to L^p(\RR)$  defined by
$$
f\mapsto D^\alpha_0 f(t):={1\over \Gamma(\alpha+1)}\vert
t\vert^\alpha W^\alpha_0f(t), \qquad t\in \R, \quad f \in
\mathcal{T}_p^{(\alpha)}(\vert t\vert^\alpha),
$$
is an isometry. Similarly we also define $D^\alpha_-$ involving $W^\alpha_-$. For
$p>1$ and $p'$ satisfies $\frac{1}{p}+\frac{1}{p'}=1$, then the dual of $\mathcal{T}_p^{(\alpha)}(\vert t\vert^\alpha)$ is $\mathcal{T}_{p'}^{(\alpha)}(\vert t\vert^\alpha)$, where the duality is given by
$$
\langle f,g \rangle_{\alpha} ={1\over
\Gamma(\alpha+1)^2}\int_{-\infty}^\infty W^\alpha_0 f(t)W^\alpha_0
g(t)\vert t\vert^{2\alpha}dt,
$$
for $f\in \mathcal{T}_p^{(\alpha)}(\vert t\vert^\alpha)$, $g\in
\mathcal{T}_{p'}^{(\alpha)}(\vert t\vert^\alpha)$, see more details in \cite[Proposition 4.2]{[LMPS]}.

As a matter of a fact, we have that $ \mathcal{T}_p^{(\alpha)}(\vert t\vert^\alpha) = \mathcal{T}_{-,p}^{(\alpha)}((-t)^\alpha) \oplus \mathcal{T}_p^{(\alpha)}(t^\alpha)$ where $\mathcal{T}_{-,p}^{(\alpha)}((-t)^\alpha)$ is a copy of $\mathcal{T}_p^{(\alpha)}(t^\alpha)$ supported on $(-\infty, 0)$, see the case $p=1$ in \cite[Section 2]{GMR}.

Now we focus our attention on the Hilbert transform on $\R$, $\mathcal{H}$. As we mention in the Introduction, it is given by
$$
{\mathcal H}f(t):=\hbox{p.v.}{i\over \pi}\int_{-\infty}^\infty{f(s)\over t-s}ds, \qquad t\in \R,
$$
and it defines a bounded operator on $L^p(\R)$ for $1<p<\infty$. Moreover, it is an isometry for $p=2$. Before giving our next result, we need the following useful Lemma.
\begin{lemma} \label{momento}For $f$ in  the Schwartz class on
$\mathbb{R}$  and $n\in \NN$, we have
$$\int_\R t^{j}f^{(n)}(t)dt=0, \qquad j\in\{0,1,\cdots n-1\}.
$$
\end{lemma}
\bgproof Note that for $0\le j\le n-1$, we have that
$$
\int_0^\infty t^{j}f^{(n)}(t)dt=\int_0^\infty (t-0)^{j}(f^{(n-j-1)})^{(j+1)}(t)dt=j!(-1)^{j+1}f^{(n-j-1)}(0),
$$
and similarly
$$
\int_{-\infty}^0 t^{j}f^{(n)}(t)dt=(-1)^j\int_{-\infty}^0 (0-t)^{j}(f^{(n-j-1)})^{(j+1)}(t)dt=j!(-1)^{j}f^{(n-j-1)}(0),
$$
and we conclude the equality.
\edproof

 Then, as
\begin{eqnarray*}
{\mathcal H}(f^{(n)})(y)&= &({\mathcal H}f)^{(n)}(y), \\
{\mathcal H}(x^nf)(y)&= &y^n{\mathcal H}f(y)-{i\over \pi}\sum_{j=0}^{n-1}y^j\int_\R t^{n-1-j}f(t)dt,
\end{eqnarray*}
for $n\in \NN$ and $y\in \RR$ (\cite[Formula (2.2)]{[Du]}), we  apply Lemma \ref{momento} to get
$$
|||{\mathcal H}f|||_{n, 2}=|||y^n({\mathcal H}f)^{(n)}|||_{0, 2}=|||{\mathcal H}(x^nf^{(n)})|||_{0, 2}= |||x^nf^{(n)}|||_{0, 2}=|||f|||_{n, 2}.
$$
We conclude that ${\mathcal H}$ is also an isometry on $\mathcal{T}_p^{(n)}( \vert t\vert^n)$ for $1<p<\infty$. In fact, similar properties hold on  $\mathcal{T}_p^{(\alpha)}( \vert t\vert^\alpha)$ for $\alpha>0$ as next result shows.

\begin{theorem}\label{Hilberttransf}For $\alpha\ge 0,$ and $1<p<\infty$.
\begin{itemize}

\item[(i)] The Hilbert transform ${\mathcal H}$ verifies that
\begin{eqnarray*}
{\mathcal H}f(t)=\left\{
	\begin{array}{ll}
		{\mathcal H}_+f_+(t)+\displaystyle{i\over \pi}{\mathcal S}\tilde{f}(t),& \mbox{if } t>0, \\
		-{\mathcal H}_+\tilde{f}(-t)-\displaystyle{i\over \pi}{\mathcal S}f_+(-t),& \mbox{if } t<0,
	\end{array}
\right.
\end{eqnarray*}
where $f_+(t):=f(t),$ $\tilde{f}(t):= f(-t)$ for $t>0$ and $f\in \mathcal{T}_p^{(\alpha)}( \vert t\vert^\alpha)$.

\item[(ii)] The equality $D_{0}^{\alpha} \mathcal{H} = \mathcal{H}  D_{0}^{\alpha}$ holds,  and as a consequence the Hilbert transform  $\mathcal{H}$ is a bounded on $\mathcal{T}_p^{(\alpha)}(\vert t\vert^\alpha)$, and $$\|\mathcal{H}\|= \cot \frac{\pi}{2p^\ast},$$
    where $p^\ast := \max\{p,p'\}$. In particular $\mathcal{H}$ is an isometry on $\mathcal{T}_2^{(\alpha)}( \vert t\vert^\alpha)$.
\end{itemize}

\end{theorem}
\bgproof (i) Take $f\in \mathcal{T}_p^{(\alpha)}(  \vert t\vert^\alpha)$. Then $f_+, \tilde{f}\in \mathcal{T}_p^{(\alpha)}(   t^\alpha)$ where $f_+(t):=f(t),$ $\tilde{f}(t):= f(-t)$ for $t>0$. Note that
$$
{\mathcal H}f(t)=\hbox{p.v.}{i\over \pi}\int_{0}^\infty{f_+(s)\over t-s}ds+{i\over \pi}\int_{0}^\infty{\tilde{f}(s)\over t+s}ds={\mathcal H}_+f_+(t)+\displaystyle{i\over \pi}{\mathcal S}\tilde{f}(t)
$$
for $t>0$ and similarly ${\mathcal H}f(t)=-{\mathcal H}_+\tilde{f}(-t)-\displaystyle{i\over \pi}{\mathcal S}f_+(-t)$ for $t<0.$

(ii) Take $t>0$. By the part (i), Lemma \ref{llave} and Theorem \ref{semihilbert} (ii) we have that
\begin{eqnarray*}
{\mathcal H}(D_0^\alpha(f))(t)&=&{\mathcal H}_+(D_+^\alpha(f_+))(t)+\displaystyle{i\over \pi}{\mathcal S}(\widetilde{D^\alpha_{-}(f)})(t)=
D_+^\alpha\left({\mathcal H}_+(f_+)(t)+\displaystyle{i\over \pi}{\mathcal S}(\tilde{f})\right)(t)\\&=&D_0^\alpha({\mathcal H}(f))(t).
\end{eqnarray*}
We prove a similar equality for $t<0$. The proof of the last part is similar to the proof of Theorem \ref{semihilbert} (ii).
\edproof

We remark that, as in the case of
$\mathcal{T}_p^{(\alpha)}(t^\alpha)$, it is easy to verify that
$(T_{t,p})_{t\in \mathbb{R}}$  where
\begin{equation}\label{grupo2}
T_{t,p} f(s):=e^{-\frac{t}{p}}f(e^{-t}s), \qquad f \in
\mathcal{T}_p^{(\alpha)}(\vert t\vert^\alpha),
\end{equation}
is a $C_0$-group of isometries on
$\mathcal{T}_p^{(\alpha)}(|t|^\alpha)$. Its is  infinitesimal generator $\Lambda$ is given by
$$
(\Lambda f)(s):=-sf'(s)-\frac{1}{p}f(s)
$$
with domain $D(\Lambda)= \mathcal{T}_p^{(\alpha+1)}(\vert
t\vert^{\alpha+1})$ and $\sigma(\Lambda)=i\mathbb{R}$.  The semigroups $(T_{t,p})_{t\geq 0}$ and $(T_{-t,p'})_{t\geq 0}$ are adjoint operators of each other acting on $\mathcal{T}_p^{(\alpha)}(\vert t\vert^\alpha)$ and $\mathcal{T}_{p'}^{(\alpha)}(\vert t\vert^\alpha)$ with  ${1\over p}+ {1\over p'}=1$ for $p>1, $ (\cite[Theorem 4.4]{[LMPS]}).

For $\mu>\beta > 0$ we define the generalized Stieltjes operator by
$$
\mathcal{S}_ {\beta,\mu} f(t) := \left\{ \begin{array}{ll} \displaystyle{|t|^{\mu-\beta}}\int_{-\infty}^0 {\vert s\vert^{\beta-1}\over \vert t+s \vert^{\mu}}f(s)ds, & t<0, \\
                                                    \\
                                                    B(\mu-\beta,\mu)f(0),                                                               & t=0, \\
                                                    \\
                                                    \displaystyle{t^{\mu-\beta}}\int^{\infty}_0 { s^{\beta-1}\over (t+ s)^{\mu}}f(s)ds,   & t>0, \\
                                  \end{array} \right.
$$
for $f\in {\mathcal S}.$ We are interested in the extension of
$\mathcal{S}_{\beta,\mu}$ on $\mathcal{T}_p^{(\alpha)}(\vert
t\vert^\alpha)$. Note that we may write
$$
\mathcal{S}_{\beta,\mu} f(t)= \int_0^\infty
{u^{\beta-1}\over (1+u)^{\mu}}f(tu)du , \quad \,\,t\in \R, \,\, f\in {\mathcal S}.
$$

We use this integral representation to prove the next lemma.

\begin{lemma}\label{llave2}
Take $\alpha \ge 0 $ and  $\mu, \beta \in \R$. Then $ D^\alpha_0
\mathcal{S}_{\beta,\mu}=\mathcal{S}_{\beta,\mu} D^\alpha_0,$ i.e.,
$$
D^\alpha_0(\mathcal{S}_{\beta,\mu}(f))=\mathcal{S}_{\beta,\mu}(
D^\alpha_0(f)),\qquad f\in {\mathcal S},
$$
where $D^\alpha_0f(t):=\frac{1}{\Gamma(\alpha+1)}\vert  t\vert
^\alpha W^\alpha_0 f(t)$ for $t\in \R$ and $f\in {\mathcal S}$.
\end{lemma}

\bgproof Since for $\lambda >0$, we have that $W_0^\alpha
(f_\lambda)= \lambda^\alpha (W_0^\alpha f)_\lambda$, where
$f_\lambda(t)= f(\lambda t)$ for $t\in \R,$ the proof follows
similarly to Lemma \ref{llave}. \edproof

Similar results of $\mathcal{S}_{\beta,\mu}$ on
$\mathcal{T}_p^{(\alpha)}(t^\alpha)$ hold for $\mathcal{S}_{\beta,\mu}$
on $\mathcal{T}_p^{(\alpha)}(\vert t\vert^\alpha)$. The proof of
next result is analogous to the proof of Theorems \ref{cotaTa},  \ref{spec} and \ref{dualStieltjes}.

\begin{theorem}\label{lemma22.1}
Let $\alpha\geq 0$, $\mu>\beta - \pp > 0$  and $1 \leq p<\infty$ and the generalized
Stieltjes operator $\mathcal{S}_{\beta,\mu}$ on
$\mathcal{T}_p^{(\alpha)}(\vert t\vert^\alpha)$. Then
\begin{itemize}
\item[(i)] The operator $\mathcal{S}_{\beta, \mu}$ is bounded  on $\mathcal{T}_p^{(\alpha)}(\vert t\vert^\alpha)$ and
$$
||\mathcal{S}_{\beta, \mu}||=B(\mu - \beta +1/p, \beta -1/p),
$$

\item[(ii)] If $f\in \mathcal{T}_p^{(\alpha)}(\vert t\vert^\alpha)$, then
$$
\mathcal{S}_{\beta,\mu} f(t)=\displaystyle\int_{-\infty}^\infty
\phi_{\mu -\beta +1/p, \mu}(r)T_{r,p}f(t)dr, \quad t \in \R,
$$
where the $C_{0}$-group $(T_{r, p})_{ r\in\R}$ is defined in (\ref{grupo2}) and the set of functions $\phi_{\beta,\mu}$ in (\ref{functionp}).

\item[(iii)]
$$
\sigma(\mathcal{S}_{\beta,\mu})= {\biggl\{ B\left(\beta-{1\over p}+i\xi,\mu-\beta+{1\over p}-i\xi\right)\,:\, \xi \in\R \biggr\}}\cup\{0\}.
$$

\item[(iv)] The adjoint operator of  generalized Stieltjes operator
$\mathcal{S}_{\beta,\mu}$  on
$\mathcal{T}_p^{(\alpha)}(\vert t\vert^\alpha)$ is $\mathcal{S}_{\mu-\beta+1,\mu} $  on $\mathcal{T}_{p'}^{(\alpha)}(\vert t\vert^\alpha)$. In particular, the operator $\mathcal{S}_{\beta,2\beta-1}$ is  self-adjoint on $\mathcal{T}_{2}^{(\alpha)}(\vert t\vert^\alpha)$ for $\beta>{1\over 2}.$
\end{itemize}
\end{theorem}

Moreover, same reasoning used to show that $D_{0}^{\alpha} \mathcal{H} = \mathcal{H}  D_{0}^{\alpha}$ in Theorem \ref{Hilberttransf} (ii), combined with the fact that $\mathcal{H}_+ \Sbmu = \Sbmu \mathcal{H}_+$ on $\Tpa$ (Theorem \ref{semihilbert} (iii)), shows that in fact the Hilbert transform and the generalized Stieltjes operator commute as operators on $\Tpabs$, that is, $\mathcal{H} \Sbmu = \Sbmu \mathcal{H}$.

We now turn to study the Fourier transform on $\Tpabs$, as well as its composition with the generalized Stieltjes operator. We first remind the reader that the Fourier transform of a function $f$
in $L^1(\RR)$ is defined by
$$
\hat f(t):=\int_{-\infty}^{\infty}e^{-ixt} f(x)dx, \qquad t\in \R.
$$
It is well-known that $\hat f$ is continuous on $\RR$ and $\hat
f(t)\to 0$ when $\vert t\vert\to \infty$ (the Riemann-Lebesgue
lemma). In the case that $f\in L^p(\RR)$ for some $1<p\le 2$, the
Fourier transform of $f$ is defined in terms of  a limit in the
norm of $L^{p'}(\RR)$ of truncated integrals:
$$
\hat f:= \lim_{R\to \infty}\widehat{f\chi_{(-R, R)}}, \qquad
\widehat{f\chi_{(-R, R)}}(t)=\int_{-R}^{R}e^{-ixt} f(x)dx, \qquad
t\in \R,
$$
i.e., $\hat f \in L^{p'}(\RR)$ and $ \lim_{R\to \infty}\Vert \hat
f-\widehat{f\chi_{(-R, R)}}\Vert_{p'}=0 $ where ${1\over
p}+{1\over p'}=1$ and $\chi_{(-R, R)}$ is the characteristic
function of the interval $(-R, R)$, see for example  \cite[Vol 2,
p.254]{[Zi]}. Then  the existence of $\hat f(t)$ is guaranteed
only at almost every $t$ and $\hat f$ may be non continuous and
the Riemann-Lebesgue lemma  could not hold (unlike the case when
$f\in L^1(\R))$.

In case that $f\in L^p(\R)$ for some $2<p<\infty$, the  Fourier
transform $\hat f$ cannot be defined as an ordinary function
although $\hat f$ can be defined as a tempered distribution, see
for example \cite[pp 19-30]{[SW]}.

We also may consider the Fourier transform on the
Sobolev spaces $\mathcal{T}_{p}^{(n)}(\vert t\vert^n)$, where $n$ is a natural number. Take $f\in
\mathcal{T}_{p}^{(n)}(\vert t\vert^n)$, then it is known that $\hat f \in \mathcal{T}_{p'}^{(n)}(\vert t\vert^n)$ for $1\le p\le 2$, $n\in \NN$ and ${1\over p}+{1\over
p'}=1$ (\cite[Theorem 6.1]{[LMPS]}). In the same work it is shown the following interesting property, about the $C_0$-group of operators  $(T_{t,p})_{t\in
\RR}$ defined by (\ref{grupo2}): given
$ f\in\mathcal{T}_p^{(\alpha)}(\vert
t\vert^\alpha),$  then
\begin{equation}\label{fg}
\widehat{T_{t,p} (f)}=T_{-t,p'}( \hat f),
\qquad \alpha \ge 0, \qquad {1\over p}+{1\over p'}=1
\end{equation}

In the next theorem, we show that the generalized Sieltjes operator and Fourier transform commute.  Our proof is based in the integral
representations of ${\mathcal S}_{\beta,\mu}(f)$  given in Theorem \ref{lemma22.1} (ii).

\begin{theorem}\label{conmutan2} For $1\le p\le 2$ and  $0< \beta -1/p < \mu$, the following equality holds
$$
\widehat{{\mathcal S}_{\beta, \mu}(f)}={\mathcal S}_{\mu-\beta+1,\mu}(\widehat{f}), \qquad
f\in L^p(\R),
$$
in particular $\widehat{{\mathcal S}_{\beta,\mu}(f)}(x)={\mathcal S}_{\mu-\beta+1,\mu}(\widehat{f})(x)$ for almost every $x$ on $\R$. In the case $n\in \NN$, $\widehat{{\mathcal S}_{\beta, \mu}(f)}={\mathcal S}_{\mu-\beta+1,\mu}(\widehat{f})$ for $f \in \mathcal{T}_p^{(n)}(\vert t\vert^n)$.
\end{theorem}

\bgproof Take $f\in L^p(\R)$ for  $1\le p\le 2$. By Theorem
\ref{lemma22.1} (ii) and formula (\ref{fg}) we have that
\begin{eqnarray*}
\widehat{{\mathcal S}_{\beta,\mu}(f)}(x)&=&
\int_{-\infty}^\infty
\phi_{\mu -\beta +1/p, \mu}(r)\widehat{T_{r,p}f}(x)dr=\int_{-\infty}^\infty
\phi_{\mu -\beta +1/p, \mu}(-r){T_{r,p'}\widehat{f}(x)dr} \\
&=&\int_{-\infty}^\infty
\phi_{\beta +1/p'-1, \mu}(r){T_{r,p'}\widehat{f}(x)dr}=\int_{-\infty}^\infty
\phi_{\mu-(\mu-\beta+1) +1/p', \mu}(r){T_{r,p'}\widehat{f}(x)dr}\\&=&{\mathcal S}_{\mu-\beta+1,\mu}(\widehat{f})(x),
\end{eqnarray*}
for a.e. $x$ in $\R$, where we used that $
\phi_{\beta, \mu}(-r)=
\phi_{\mu-\beta, \mu}(r)
$ for $r\in \R$. Finally, given $f \in \mathcal{T}_p^{(n)}(\vert
t\vert^n)$, then $\widehat{f}\in \mathcal{T}_{p'}^{(n)}(\vert
t\vert^n)$, and we conclude that $\widehat{{\mathcal S}_{\beta, \mu}(f)}={\mathcal S}_{\mu-\beta+1,\mu}(\widehat{f})$.
\edproof

\section{Stieltjes convolution product}

\setcounter{theorem}{0} \setcounter{equation}{0}

In \cite{[Ya]} (see also \cite{[YaM]}) authors introduced  an interesting  convolution operator
\begin{eqnarray}\label{convol}
(f \otimes g) (t) &:=& f(t)  \int_0^\infty \frac{g(s)}{s-t}ds + g(t) \int_0^\infty \frac{f(s)}{s-t}ds, \qquad t>0,
\end{eqnarray}
where both integrals $\int_0^\infty$ are indeed improper integrals $\hbox{p.v.} \int_0^\infty$.
The following convolution-type identity for  Stieltjes
transform,
\begin{eqnarray}\label{product}
\mathcal{S}(f\otimes g) = \mathcal{S}f \cdot \mathcal{S}g,
\end{eqnarray}
holds  in a class of functions, which is associated with the Mellin transform; $(\cdot)$ stands for pointwise multiplication, and $\mathcal{S}$ for the Stieltjes transform.

 Later these results were
extended on $L^p$-spaces, $f \otimes g\in L^r(\R^+)$ in the case that $f\in L^p(\R^+)$ and $ g \in L^q(\R^+)$ when ${{1\over p}+{1 \over q} = {1\over r}}$ (\cite{[SV]}). Then  they are applied to a class of singular integral equations of convolution type.

Note that the equality (\ref{convol}) may be written as
$$
f \otimes g  :=f{\mathcal H}_+g + g{\mathcal H}_+f,
$$
where the operator ${\mathcal H}_+$ is defined in Theorem \ref{Hilberttransf} (i). Recall that $\mathcal{S} = \mathcal{S}_{1,1}$, and a natural question arises: does the generalized Stieljes operator $\mathcal{S}_{\beta, \mu}$ (defined on Lebesgue spaces $L^p(\R^+)$ and more generally in Sobolev-Lebesgue spaces, $\mathcal{T}_p^{(\alpha)}(\vert
t\vert^\alpha),$ for $ p>1$) satisfy a similar convolution-type identities to (\ref{product})  for arbitrary positive $\beta, \mu$?  We give a partial answer (for  $\beta, \mu$ natural numbers) on the next theorem, but first we will need the Lemma below:

\begin{lemma}\label{general} Let $p, q  \in (1, \infty)$  and $r\ge 1$ be such that ${1 \over p}+{1 \over q} = {1 \over r},$
  $f \in L^p(\R^+)$ and $g \in L^q(\R^+)$. Then
  \begin{eqnarray*}
  {\mathcal S}_{\beta, \mu}(g{\mathcal H}_+f)(t)&=& t^{\mu-\beta}\int_0^\infty f(u)\left(\hbox{\normalfont{p.v.}}\int_0^\infty \frac{s^{\beta-1}}{(s+t)^\mu(u-s)} g(s)ds\right)du,\\
\mathcal{S}_{\beta,\mu}(f\otimes g) (t) &=&   t^{\mu-\beta} \int_0^\infty {f(s)\over (t+s)^\mu}\left(\hbox{\normalfont{p.v.}}\int_0^\infty h_{\beta, \mu}(t,s,u){g(u)\over (t+u)^\mu}du\right)ds,
\end{eqnarray*}
where $h_{\beta, \mu}(t,s,u)=\displaystyle{s^{\beta-1}(t+u)^\mu- u^{\beta-1}(t+s)^\mu\over {u-s}}$ for $0 <\beta- {1\over r} <\mu$ and $t>0$.
\end{lemma}

\bgproof  By the classical H$\ddot{\hbox{o}}$lder inequality, given $f\in L^p(\R^+)$ and $g \in L^q(\R^+)$, then $f{\mathcal H}_+g\in L^r(\R^+)$. Then ${\mathcal S}_{\beta, \mu}(g{\mathcal H}_+f)\in L^r(\R^+)$ for $0 <\beta- {1\over r} <\mu$.

Now take $\varepsilon >0$ and  consider the function
\begin{eqnarray*}
{\mathcal H}_{+,\varepsilon}f(s) :=\int_{(0,s-\varepsilon)\cup(s+\varepsilon,\infty)} \frac{f(u)}{u-s}du, \qquad s>0.
\end{eqnarray*}
Note that   ${\mathcal H}_{+,\varepsilon}f \in L^p(\R^+)$, ${\mathcal H}_{+,\varepsilon}f(s)\to{\mathcal H}_+ f(s)$ for $s$ a.e. and $\Vert {\mathcal H}_{+,\varepsilon}f-{\mathcal H}_+ f\Vert_p\to 0$ when $\varepsilon \to 0$  (\cite[Theorem 8.1.12]{BN}). By the boundedness of the Hilbert transform, we conclude that
$$
\Vert {\mathcal H}_{+,\varepsilon}f\Vert_p\le \Vert   {\mathcal H}_{+,\varepsilon}f-{\mathcal H}_+ f\Vert_p+\Vert{\mathcal H}_+ f \Vert_p\le C \Vert f\Vert_p.
$$

Fixed $t>0$, we consider the function $\kappa_t(s):=t^{\mu-\beta}\displaystyle{\frac{s^{\beta-1}}{(t+s)^\mu}} \in L^\rho(\R^+)$ for ${1 \over r} + {1 \over \rho} = 1$. Now we apply the H\"{o}lder inequality for three functions $\kappa_t,$ $g$ and ${\mathcal H}_{+,\varepsilon}f$ to conclude that  $\omega_{\varepsilon,t} :=  \kappa_tg{\mathcal H}_{+,\varepsilon}f\in L^1(\R^+)$ and
$$
\Vert \omega_{\varepsilon,t}\Vert_1\le \Vert \kappa_t\Vert_\rho\Vert g\Vert_q\Vert {\mathcal H}_{+,\varepsilon}f\Vert_p\le C\Vert \kappa_t\Vert_\rho\Vert g\Vert_q\Vert f\Vert_p.
$$
On the other hand $\vert \omega_{\varepsilon,t}(s)\vert \le C \kappa_t(s) \vert g (s)\vert\,{\mathcal H}_+^\ast f(s)  $  for $s$ a.e.
where $\mathcal{H}_+^\ast$ is  the maximal operator defined as $\mathcal{H}_+^\ast f (x) = \sup_{\varepsilon>0}|\mathcal{H}_{+,\varepsilon}f(x)|$.
By Lebesgue dominated convergence theorem, we have that
\begin{eqnarray*}
  {\mathcal S}_{\beta, \mu}(g{\mathcal H}_+f)(t)&=&
 \int_0^\infty \lim_{\varepsilon \rightarrow 0^+} \omega_{\varepsilon,t}(s)ds
 = \lim_{\varepsilon \rightarrow 0^+} \int_0^\infty \omega_{\varepsilon,t}(s)ds\\&=&t^{\mu-\beta}\lim_{\varepsilon \rightarrow 0^+} \int_0^\infty \frac{s^{\beta-1}}{(s+t)^\mu} g(s){\mathcal H}_{+,\varepsilon}f(s)ds.
\end{eqnarray*}
By Fubini theorem, and putting the improper integrals $\left(\int_0^{t-\varepsilon} +\int_{t+\varepsilon}^\infty\right)$ simply as $\int_0^\infty$
\begin{eqnarray*}
t^{\mu-\beta} \int_0^\infty \frac{s^{\beta-1}}{(s+t)^\mu} g(s)\int_0^\infty \frac{f(u)}{u-s}duds
= t^{\mu-\beta}\int_0^\infty f(u)\left(\int_0^\infty \frac{s^{\beta-1}}{(s+t)^\mu(u-s)} g(s)ds\right)du.
\end{eqnarray*}

Now we use  this identity to get
\begin{eqnarray*}
\mathcal{S}_{\beta,\mu}(f\otimes g) (t) &=& t^{\mu-\beta}\int_0^\infty \frac{s^{\beta-1}}{(s+t)^\mu}\left(
f(s) \int_0^\infty \frac{g(u)}{u-s}du + g(s) \int_0^\infty \frac{f(u)}{u-s}du \right)
\\
&=&  t^{\mu-\beta} \int_0^\infty f(s) \int_0^\infty \frac{g(u)}{u-s}  \left( \frac{s^{\beta-1}}{(t+s)^\mu} -\frac{u^{\beta-1}}{(t+u)^\mu}\right)duds\\
&=&   t^{\mu-\beta} \int_0^\infty {f(s)\over (t+s)^\mu}\left(\int_0^\infty h_{\beta, \mu}(t,s,u){g(u)\over (t+u)^\mu}du\right)ds,
\end{eqnarray*}
where $h_{\beta, \mu}(t,s,u)=\displaystyle{s^{\beta-1}(t+u)^\mu- u^{\beta-1}(t+s)^\mu\over {u-s}}$ and we obtain the result.
\edproof

\begin{theorem}\label{convProd}
Let $p, q \in (1, \infty), \, r \geq 1$ be such that ${1 \over p}+{1 \over q} = {1 \over r}$ and $n,m \in \NN$ be such that $0 < n -\frac{1}{r} < m$, and  Then, for any $f \in L^p(\R^+)$, $g \in L^q(\R^+)$, we have that
$$
\mathcal{S}_{n,m} (f\otimes g) = \sum_{i=n}^m {m \choose i} \sum_{j=0}^{i-n}  \mathcal{S}_{i-j,m}f \mathcal{S}_{n+j,m}g-\sum_{i=0}^{n-2} {m \choose i}\sum_{j=0}^{n-2-i}  \mathcal{S}_{n-j-1,m}f \mathcal{S}_{i+j+1,m}g
$$
where $f\otimes g \in L^r(\R^+)$.
\end{theorem}

\bgproof As   $m,n \in \NN$, we apply the  Newton binomial and cyclotomic formula in the numerator of the rational function $h_{n, m}$ (Lemma {\ref{general}}) to get
\begin{eqnarray*}
&&s^{n-1}(t+s)^m - u^{n-1}(t+u)^m = s^{n-1}\sum_{i=0}^m {m \choose i} t^{m-i} u^{i} - u^{n-1}\sum_{i=0}^m {m \choose i} t^{m-i} s^{i}
\\
&=& \sum_{i=0}^m\hbox{sign}(i-n+1){m \choose i}t^{m-i}(us)^{\min\left\lbrace i,n-1\right\rbrace} (u^{|i-n+1|}-s^{|i-n+1|})
\\
&=& (u-s) \sum_{i=0, i \neq n-1}^m \hbox{sign}(i-n+1){m \choose i}t^{m-i}(us)^{\min\left\lbrace i,n-1\right\rbrace} \sum_{j=0}^{|i-n+1|-1} u^{j}s^{|i-n+1|-1-j}.
\end{eqnarray*}
Combining all these identities in  Lemma \ref{general}, one finally gets that
\begin{eqnarray*}
\noindent &\,&\mathcal{S}_{n,m}(f\otimes g) (t)=\qquad \qquad\\
\\ &=& - \sum_{i=0}^{n-2} \sum_{j=0}^{n-2-i} {m \choose i} \int_0^\infty \int_0^\infty t^{m-(n-j-1)}\frac{s^{n-j-2}}{(t+s)^m}f(s) t^{m-(i+j+1)}\frac{u^{i+j}}{(t+u)^m}g(u)duds
\\ &\,&\qquad + \sum_{i=n}^m \sum_{j=0}^{i-n} {m \choose i} \int_0^\infty \int_0^\infty t^{m-(i-j)}\frac{s^{i-j-1}}{(t+s)^m}f(s) t^{m-(n+j)}\frac{u^{n+j-1}}{(t+u)^m}g(u)duds\\
&=&\sum_{i=n}^m {m \choose i} \sum_{j=0}^{i-n}  \mathcal{S}_{i-j,m}(t)f \mathcal{S}_{n+j,m}g(t)-\sum_{i=0}^{n-2} {m \choose i}\sum_{j=0}^{n-2-i}  \mathcal{S}_{n-j-1,m}f \mathcal(t){S}_{i+j+1,m}g(t),
\end{eqnarray*}
and the proof is concluded. \edproof

\begin{remark}

\begin{enumerate}
	\item[(i)]{\rm An alternative expression of $\mathcal{S}_{n,m} (f\otimes g) $  is the following
$$
\mathcal{S}_{n,m} (f\otimes g) = \sum_{i=0}^{m-n} \mathcal{S}_{n+i,m}f
\sum_{j=n+i}^m {m \choose j} \mathcal{S}_{j-i,m}g
-\sum_{i=1}^{n-1}\mathcal{S}_{n-i,m}f \sum_{j=0}^{n-1-i}{m \choose j} \mathcal{S}_{i+j,m}g.
$$
Some particular expression of this general formula for $r >1$ are the following
\begin{eqnarray*}
\mathcal{S}_{1,1} (f \otimes g) & =& \mathcal{S}_{1,1}f \mathcal{S}_{1,1}g,\\
\mathcal{S}_{2,2} (f \otimes g) & =& \mathcal{S}_{2,2}f \mathcal{S}_{2,2}g- \mathcal{S}_{1,2}f \mathcal{S}_{1,2}g,\\
\mathcal{S}_{1,2} (f \otimes g) & =& \mathcal{S}_{2,2}f \mathcal{S}_{1,2}g+ \mathcal{S}_{1,2}f \mathcal{S}_{2,2}g + 2\mathcal{S}_{1,2}f \mathcal{S}_{1,2}g.
\end{eqnarray*}
By Theorem \ref{dualStieltjes} $(\mathcal{S}_{1,2})'= \mathcal{S}_{2,2}$ and $(\mathcal{S}_{2,2})' = \mathcal{S}_{1,2}$ and
\begin{eqnarray*}
	\mathcal{S}_{2,2} (f \otimes g) &=& \mathcal{S}_{2,2}f \mathcal{S}_{2,2}g
	- (\mathcal{S}_{2,2})'f (\mathcal{S}_{2,2})'g,\\
	\mathcal{S}_{1,2} (f \otimes g)  &=& (\mathcal{S}_{1,2})' f \mathcal{S}_{1,2}g
	+ \mathcal{S}_{1,2}f (\mathcal{S}_{1,2})' g + 2\mathcal{S}_{1,2}f \mathcal{S}_{1,2}g.
\end{eqnarray*}
}

\item[(ii)] {\rm Since $\Tpa \hookrightarrow L^p(\R^+)$ for $p \geq 1, \, \alpha \geq 0$, and by Theorem \ref{Holder}, classical Hölder inequality also holds in $\Tpa$ spaces $\left(\frac{1}{p}+\frac{1}{q}=1\right)$, then the identity given in Theorem \ref{convProd} above it is also valid as a closed expression on $\Tpa$ spaces when either $r=1$, or $r>1$ and $\alpha \in \NN$ (see Section 8.3 (ii)).}
\end{enumerate}
\end{remark}

\section{Spectral pictures, final comments}

The main aim of this last section is to illustrate and visualize some results which were proved in the paper. We also give some ideas, comments and open problems to continue our research.

\subsection{Spectral pictures}

We will draw  the spectrum sets of the Stieltjes operators, $\sigma(\mathcal{S}_{\beta,\mu}),$ on ${\mathcal T}_p^{(\alpha)}(t^\alpha)$ in some particular cases. We will use the software Mathematica in order to present these spectrum sets.  In \cite[Section 8]{AM} the spectrum of generalized discrete Ces\`aro operators pictures have been represented. As it is commented there, those curves  are also the spectrum of
continuous Ces\`aro operators on ${\mathcal T}_p^{(\alpha)}(t^\alpha)$.

By simplicity, for each $\gamma>0$ and $\mu>\gamma,$ we consider  the curve $$\Upsilon_{\gamma,\mu}:={\biggl\{ B(\gamma+i\xi,\mu-\gamma-i\xi)\,:\, \xi\in\R \biggr\}}.$$
As we mention in Theorem \ref{spec}, the spectrum of  the Stieltjes operator $\mathcal{S}_{\beta,\mu}$ on ${\mathcal T}_p^{(\alpha)}(t^\alpha)$ is the curve
$\Upsilon_{\beta-{1\over p},\mu}\cup\{0\}$, i.e.,
$$
\sigma(\mathcal{S}_{\beta,\mu})=\Upsilon_{\beta-{1\over p},\mu}\cup\{0\}.
$$

 For each $\gamma>0$ and $\mu>\gamma,$   the curve $\Upsilon_{\gamma,\mu}$ is symmetrical with respect to the OX axis, $\Upsilon_{\gamma,\mu}=\Upsilon_{\mu-\gamma,\mu},$ takes the point $B(\mu-\gamma, \gamma)$ on the complex plane (doing $\xi=0$). Moreover note that the curve $\Upsilon_{\gamma,\mu}$ is contained in the circle of center $(0,0)$ and radio $B(\mu-\gamma, \gamma)$, due to,
 $$
 \vert B(\gamma+i\xi,\mu-\gamma-i\xi)\vert \le B(\gamma,\mu-\gamma), \qquad \xi\in \R.
 $$
 On the other hand, $B(\gamma+i\xi,\mu-\gamma-i\xi)\to 0$ when $\xi\to\pm\infty,$ due to
\begin{equation}\label{double3}
\frac{\Gamma(z+\alpha)}{\Gamma(z)}=z^{\alpha}\left(1+O\left({1\over |z|}\right)\right), \quad z\in\C_+,\,\Re\alpha>0,
\end{equation}
 see, for example, \cite{ET}.

 Fixed $\mu>0$, note that the function $\gamma \mapsto B(\gamma,\mu-\gamma)$,  $(0,\mu)\to \R^+$ has a minimum at $\gamma={\mu\over 2}$, $$B\left({\mu\over 2},{\mu\over 2}\right)={\Gamma^2({\mu\over 2})\over \Gamma(\mu)}, $$ and $B(\gamma,\mu-\gamma)\to \infty $ when $\gamma \to 0^+, \mu^-$.

\begin{figure}[h]
\caption{}
\centering
\includegraphics[width=0.6\textwidth]{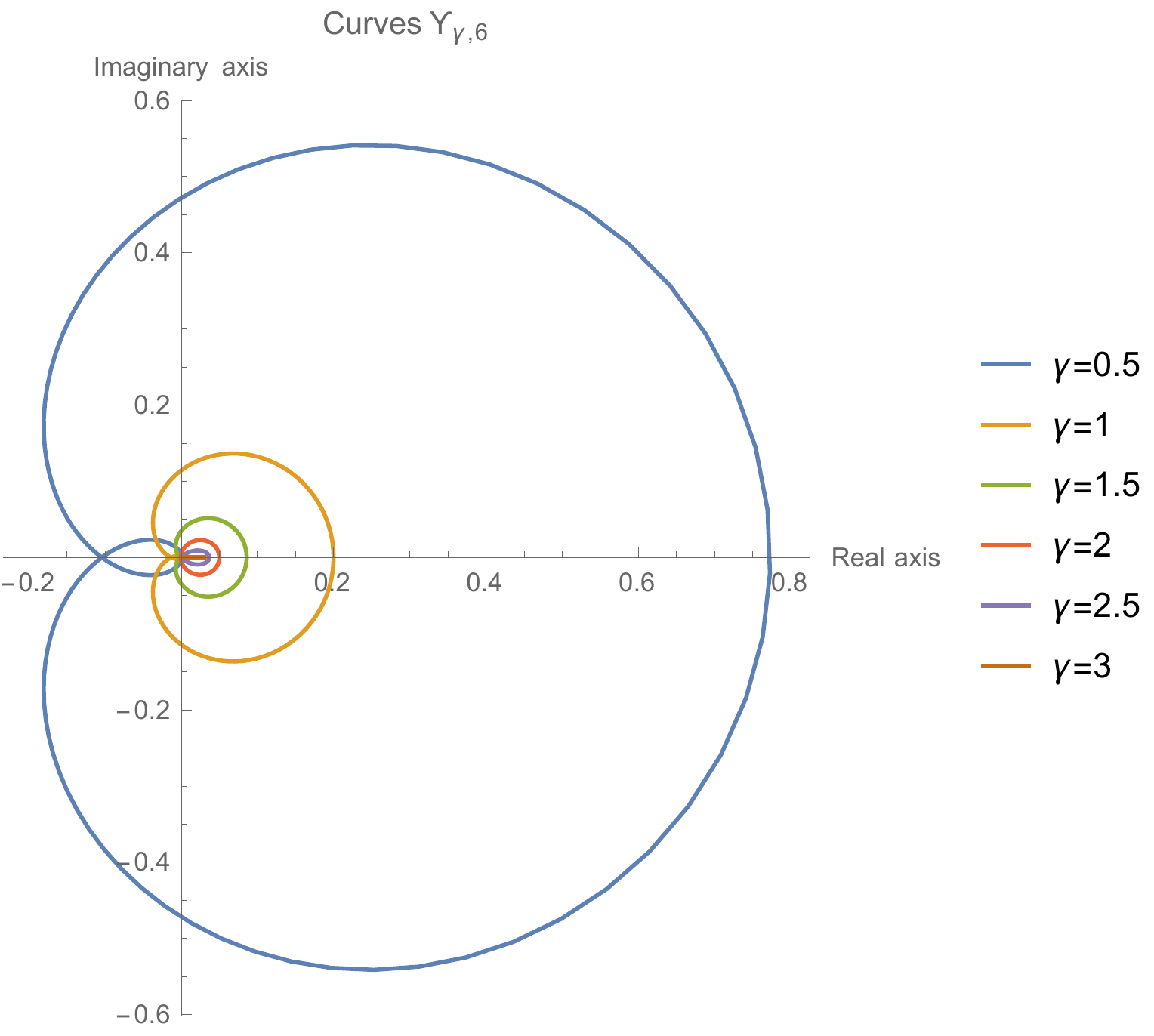}
\end{figure}

 Fixed $\gamma>0$, note that  $B(\gamma,\mu-\gamma)\to 0$ when $\mu\to \infty $ (we write $\Upsilon_{\gamma, \infty}=\{0\}$) and $B(\gamma,\mu-\gamma)\to \infty$ in the case that $\mu\to \gamma^+$.

\begin{figure}[h]
\caption{}
\centering
\includegraphics[width=0.6\textwidth]{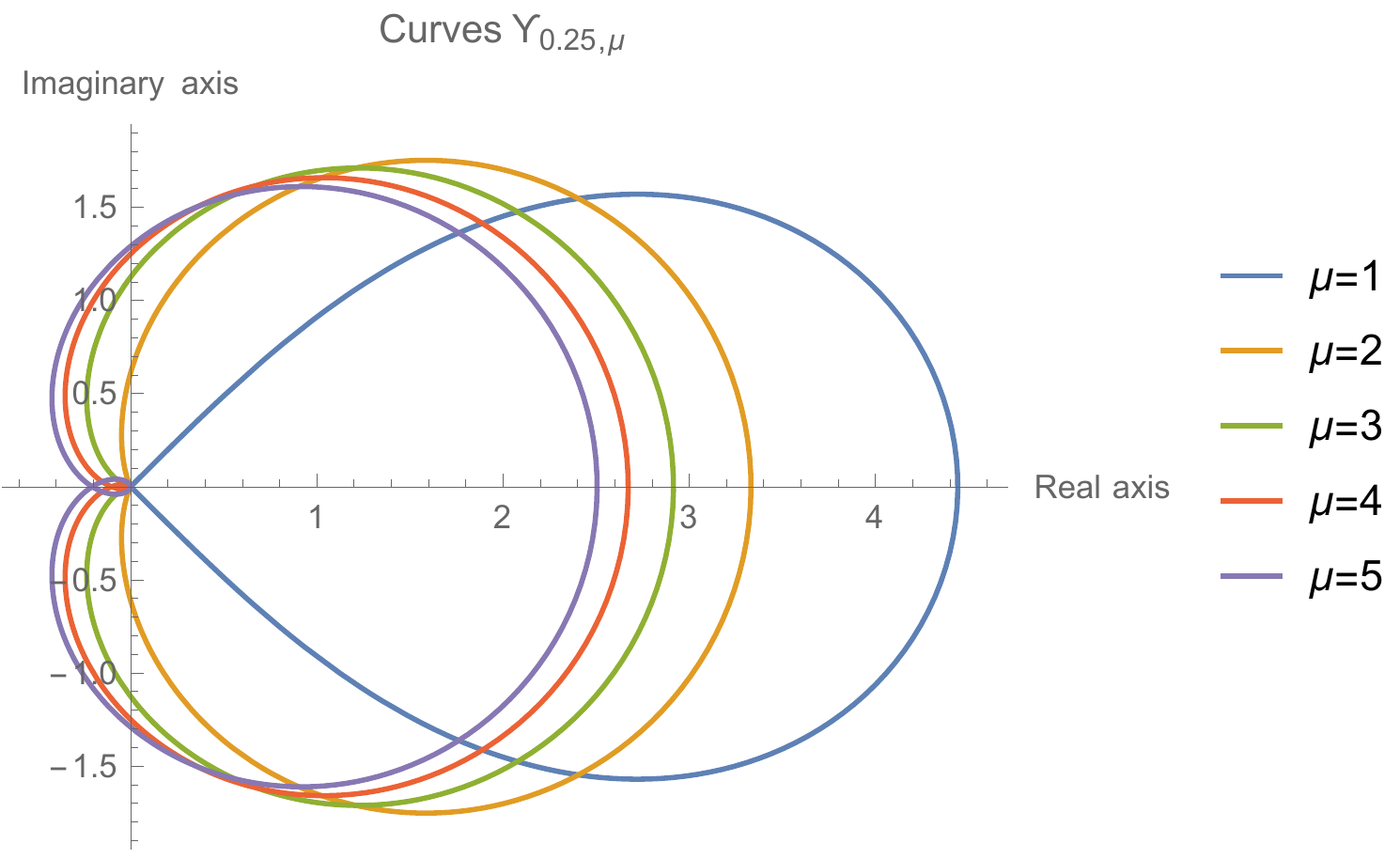}
\end{figure}

 The special case $\mu=2\gamma$ has special properties.  Since $\Gamma(\overline z)=\overline{\Gamma(z)}$, then $B(\gamma+i\xi,\gamma-i\xi)\ge 0$  for $\xi\in \R$ and $\Upsilon_{\gamma, 2\gamma}=[0, B(\gamma, \gamma)]$. From here, we may consider the extreme cases $\Upsilon_{0, 0}=[0,\infty)$ and $ \Upsilon_{\infty, \infty}=\{0\}.$  Note that, as we comment in the Remark \ref{autoad}, for $p=2$, the spectrum of $\mathcal{S}_{\beta,2\beta-1}$ is $[0, B(\beta-{1\over 2}, \beta-{1\over 2})]$, and in particular for $\beta=1$, $\sigma(\mathcal{S}_{1,1})=\sigma(\mathcal{S})=[0,\pi]$.

  However this is not the only case that the spectrum set $\sigma(\mathcal{S}_{\beta,\mu})$ is a bounded positive interval on $\mathcal{T}_p^{(\alpha)}(t^\alpha)$ for $\alpha>0$ with $p\not=2$. Take $1<p<\infty$, $ \beta>{1\over p}$ and $\mu= 2\beta-{2\over p}$. By Theorem \ref{spec}, the spectrum set of operator $\mathcal{S}_{\beta,2\beta-{2\over p}}$ is equal $[0, B(\beta-{1\over p}, \beta-{1\over p})]$  on the space $\mathcal{T}_p^{(\alpha)}(t^\alpha)$ for $\alpha\ge 0$, in particular, on $L^p(\R^+)$.

For $\mu=1$, and $0<\gamma<1$ note that
\begin{eqnarray*}
\Gamma(\gamma+i\xi)\Gamma(1-\gamma-i\xi)&=&{\pi\over \sin(\pi(\gamma+i\xi))}\cr
&=&\pi \frac{\sin(\pi\gamma)\cosh(\pi \xi)-i\cos(\pi \gamma)\sinh(\pi \xi)}{\sin^2(\pi \gamma) \cosh^2(\pi \xi)+ \cos^2(\pi\gamma)\sinh^2(\pi \xi)},
\end{eqnarray*}
and we conclude that $\Upsilon_{\gamma, 1}\subset \C^+$, (Figure 3). For $\gamma={1\over 4}$, we have that
$$
\Upsilon_{{1\over 4}, 1}={\biggl\{{\sqrt{2}\pi\over \cosh(2\pi \xi)}\left(\cosh(\pi \xi)-i\sinh(\pi \xi)\right)\,:\, \xi\in\R \biggr\}}.
$$

\begin{figure}[h]
\caption{}
\centering
\includegraphics[width=0.6\textwidth]{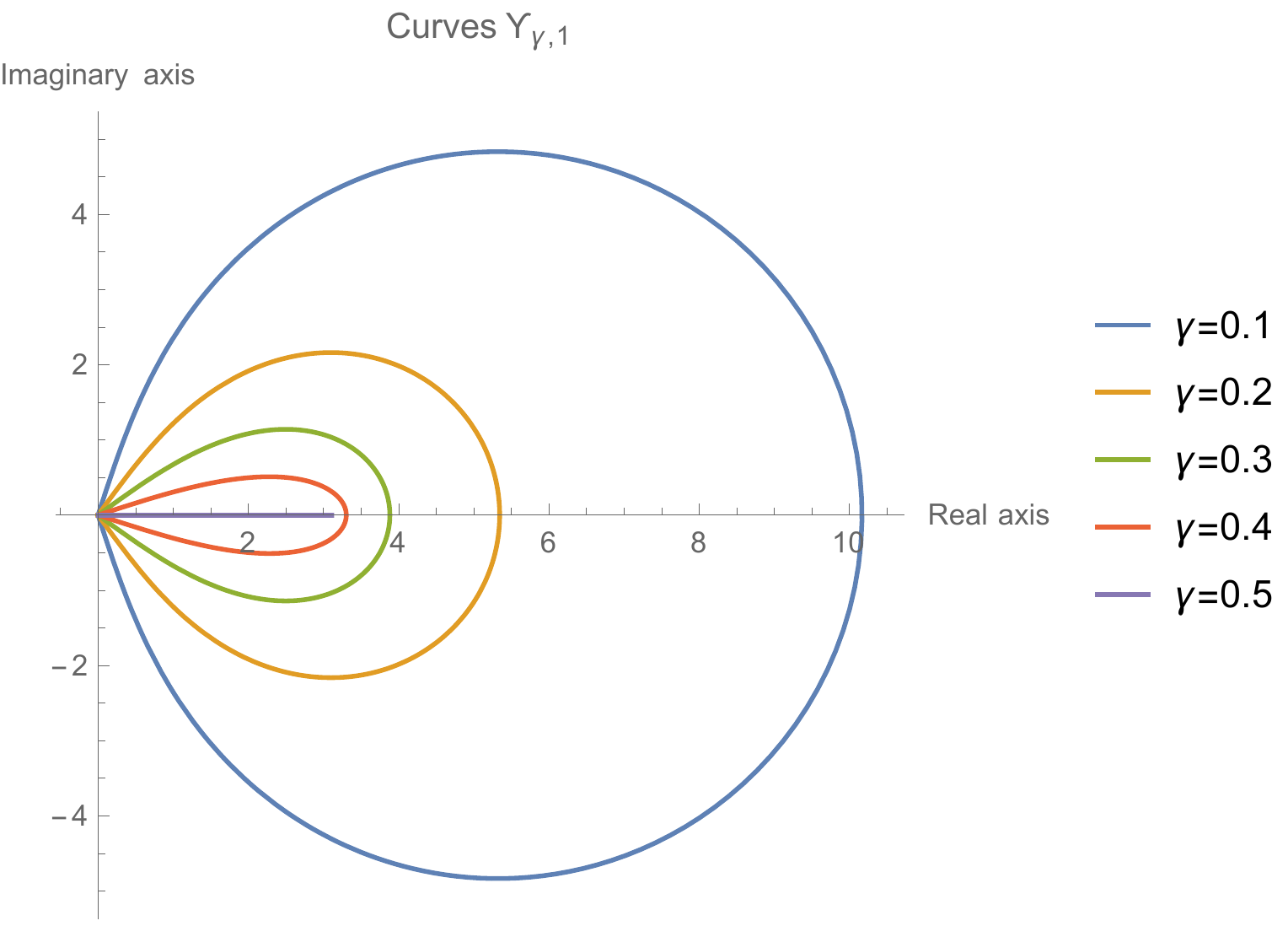}
\end{figure}

For $\mu=2$, and $0<\gamma<2$ note that
$$
\Gamma(\gamma+i\xi)\Gamma(2-\gamma-i\xi)=\pi (1-\gamma-i\xi)\frac{\sin(\pi\gamma)\cosh(\pi \xi)-i\cos(\pi \gamma)\sinh(\pi \xi)}{\sin^2(\pi \gamma) \cosh^2(\pi \xi)+ \cos^2(\pi\gamma)\sinh^2(\pi \xi)},
$$
and we do not conclude that  $\Upsilon_{\gamma, 2}\subset \C^+$. In particular for $\gamma={1\over 4}$, we obtain that,
$$
\Upsilon_{{1\over 4}, 2}={\biggl\{{\sqrt{2}\pi\over \cosh(2\pi \xi)}\left({3\over 4}\cosh(\pi \xi)-\xi\sinh(\pi \xi)-i\left(\xi\cosh(\pi \xi)+{3\over 4}\sinh(\pi \xi)\right)\right)\,:\, \xi\in\R \biggr\}},
$$
and  $\Upsilon_{{1\over 4}, 2}\not \subset \C^+$, see Figure 4.

\begin{figure}[h]
\caption{}
\centering
\includegraphics[width=0.6\textwidth]{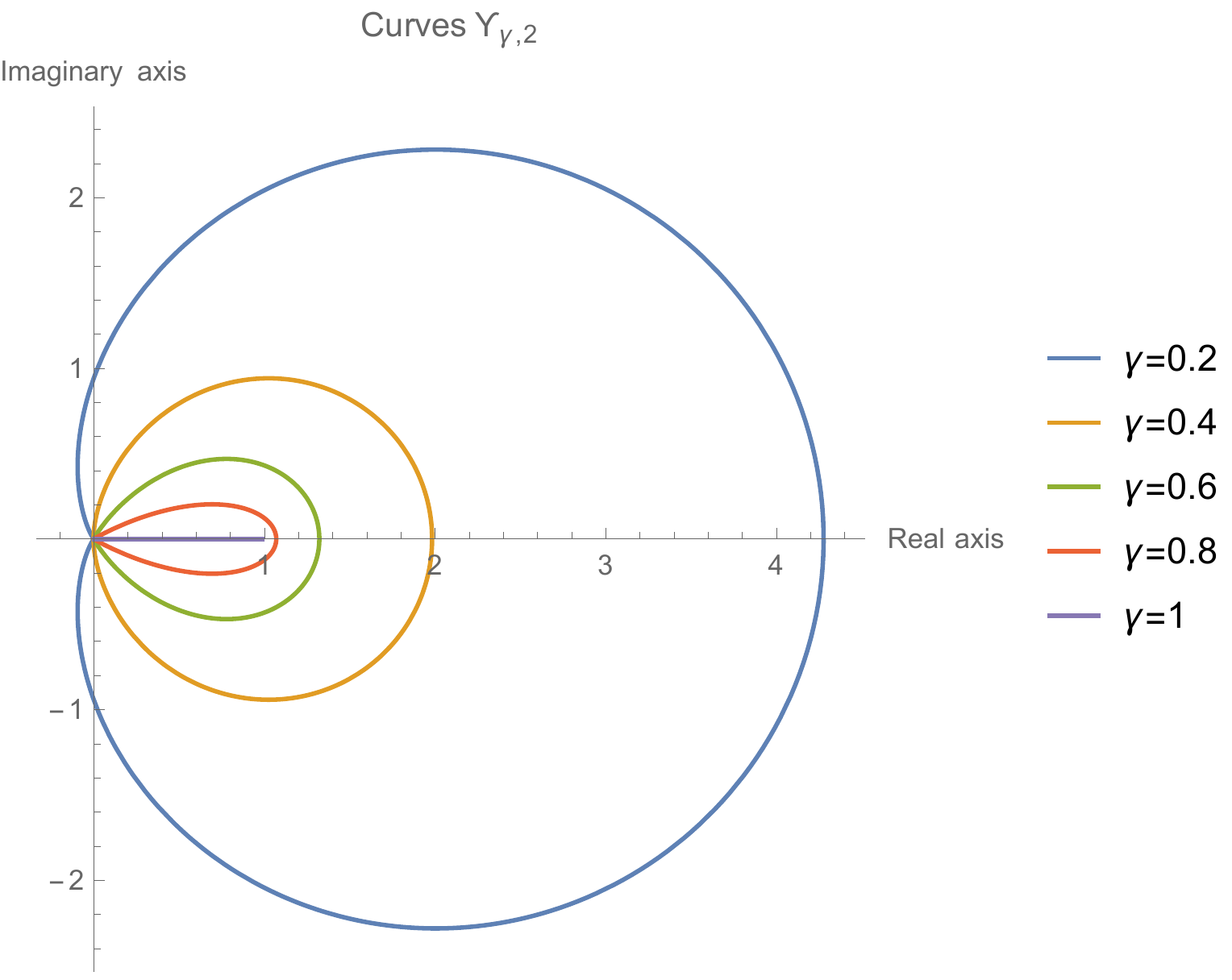}
\end{figure}

When $\mu \to \infty$, the curve $\Upsilon_{1, \mu}$ cuts to the real axis several times. Except the first cut $(t=0)$,  at the point $B(1,\mu-1),$ every cut to the real axis has double multiplicity, see Figure 5. Note that the curve $\Upsilon_{1, \mu}$ collapses into $\{0\}$ when $\mu \to \infty$.

\begin{figure}[h]
\caption{}
\centering
\includegraphics[width=0.5\textwidth]{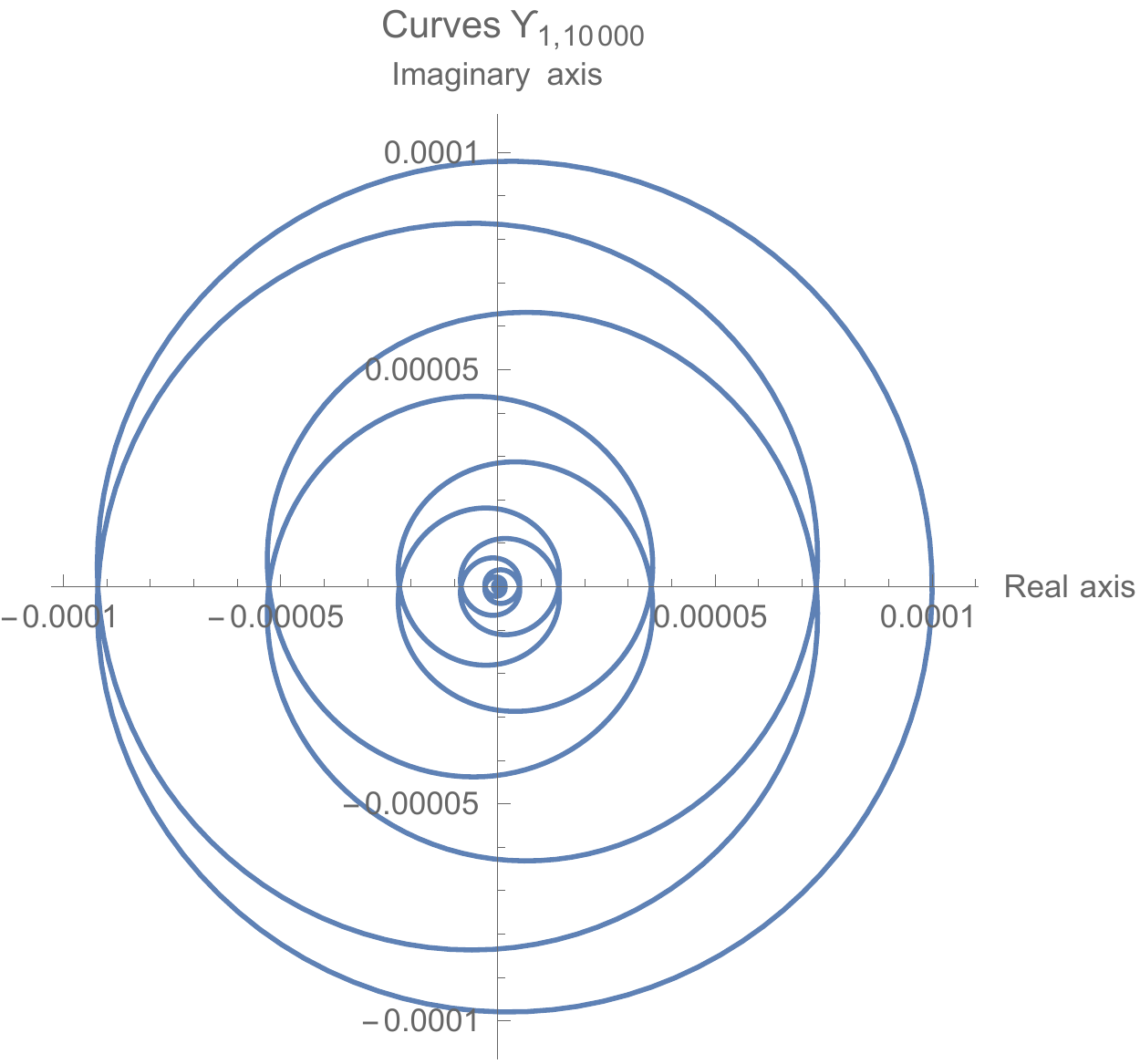}
\end{figure}

Finally we may consider the curves $\Upsilon_{\gamma, \mu}$ with $0<\Re\gamma<\Re \mu$. Note that $\Upsilon_{\gamma, \mu}= \Upsilon_{\Re\gamma, \mu}$. However it is more difficult to describe $\Upsilon_{\gamma, \mu}$ in terms of $\Upsilon_{\gamma, \Re\mu}$, compare $\Upsilon_{\gamma, 1}$ and $\Upsilon_{\gamma, 1+ 2\pi i}$ (Figures 1 and 6).

\begin{figure}[h]
\caption{}
\centering
\includegraphics[width=0.6\textwidth]{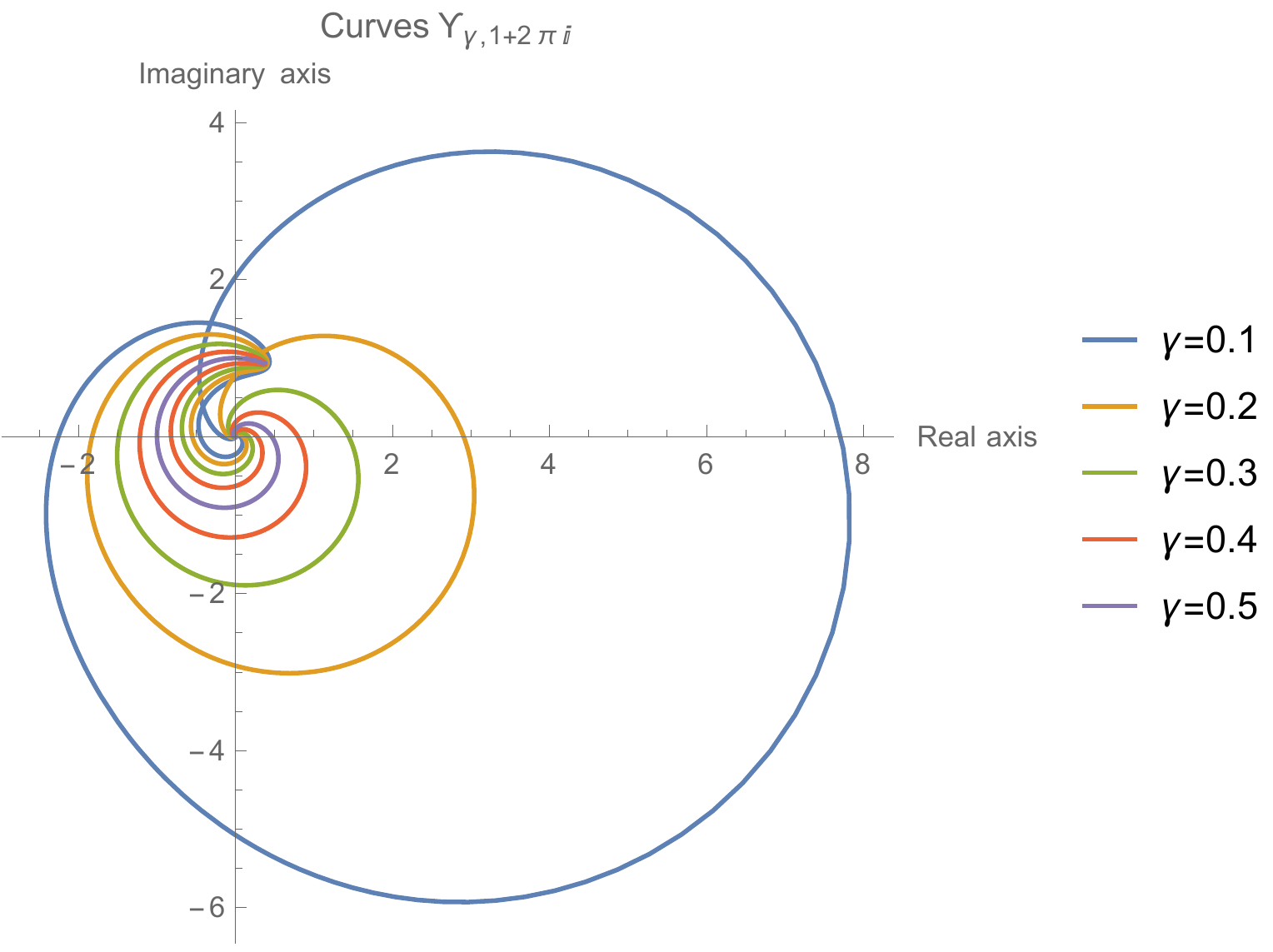}
\end{figure}
\subsection{Stieltjes operator of analytic functions}

Let $H_p$ denote the space of analytic functions $F:\C^+\to \C$  satisfying condition
\begin{equation}\label{Norm01}
\| F\| _p:=\sup_{x>0} \left(\frac{1}{\pi }\int_{-\infty }^{\infty }|F(x+iy)|^p\, dy\right)^{1\over p} <\infty .
\end{equation}
For $p=\infty$, we denote by $H_\infty$ the space of all bounded analytic functions on $\C^+$ with the supremum norm. The spaces $H_p$ , $1\le p\le \infty,$ are Banach spaces and for $p=2$, $H_2$ is a Hilbert space, used called the Hardy space on the half plane, see more details in \cite[Chapter 8]{Hoffman}.

The classical Paley-Wiener theorem states that the Laplace transform ${\mathcal L}: L^2(\R^+) \to H_2$,
$$
{\mathcal L}(f)(z)=\int_0^\infty f(t)e^{-zt}dt, \quad f\in L^2(\R^+), \quad z\in \C^+,
$$
is  an isometric isomorphism; i.e., $F, G\in H_2 $ if and only if there exist unique $ f,g\in L^2(\R^+)$ such that
$ F=\LL f$ and  $ G=\LL g$ and
$
\left\langle f,g \right\rangle=\langle F, G \rangle$, see for example \cite[Theorem 9.13]{Rud87}.

The Stietjes operator ${\mathcal S}$ (called the continuous Hilbert operator in \cite{[ASV]}) given by
$$
{\mathcal S}(F)(z):=\int_0^\infty{F(s)\over s+z}ds,\qquad z\in \C^+,
$$
is a bounded operator on $H_p$, $1<p<\infty$,  $\Vert {\mathcal S}\Vert =\displaystyle{\pi \over \sin({\pi \over p} )}$ and
$$
\sigma({\mathcal S})=\overline{\left\{{\pi\over \sin\left({\pi \over p}+i\pi \xi\right)}\,\,:\,\,\xi\in \RR\right\}},
$$
see \cite{[ASV]}.

In  \cite{[GMMS]}, fixed $n \in \NN$, the subspaces $H^{(n)}_2$ which are formed by  all analytic functions $F$ on $\C^+$ such that
\begin{equation}
z^kF^{(k)}\in H_2, \quad k=0,1,\dots,n,
\end{equation}
are introduced. These spaces are of reproducing kernels, and are  obtained as ranges of the Laplace transform in extended versions of the Paley-Wiener theorem, ${\mathcal L}({\mathcal T}_2^{(n)}(t^n))=H^{(n)}_2$ (\cite[Theorem 3.3]{[GMMS]}).

It seems natural to introduce the subspaces $H^{(n)}_p$ for $1\le p\le \infty$, and consider the generalized Stieltjes operador ${\mathcal S}_{\beta,\mu}$ given by

$$
\mathcal{S}_{\beta, \mu} F(z) := z^{\mu-\beta}\int_0^\infty \frac{s^{\beta -1}}{(z+s)^\mu}F(s)ds, \qquad z\in\C^+,
$$
on the spaces $H^{(n)}_p$ for $0< \beta -1/p < \mu$. This will be the focus of a forthcoming paper.

\subsection{ Open questions }

 In this paper we have presented a complete study of generalized Stieltjes operator on fractional Sobolev-Lebesgue spaces. However some of these results might be improved.

\begin{itemize}
\item[(i)] The case $p=\infty$ in the Sobolev-Lebesgue spaces $\mathcal{T}_{\infty}^{(\alpha)}(t^\alpha)$   is commented in Remark \ref{infinito}. It seems natural to conjecture that the operator ${\mathcal S}_{\beta, \mu}$ is bounded on these spaces $\mathcal{T}_{\infty}^{(\alpha)}(t^\alpha)$,
    $
    \Vert {\mathcal S}_{\beta, \mu}\Vert=B(\beta, \mu-\beta),$ for $0<\beta<\mu$ and $\sigma( {\mathcal S}_{\beta, \mu})=\Upsilon_{\beta,\mu}\cup\{0\}.$

\item[(ii)] In Theorem \ref{Holder} we have proved a H\"{o}lder inequality in  $\mathcal{T}_p^{(\alpha)}(t^\alpha)$ for conjugate exponents $p,p'>1$ and $\alpha \ge 0$. As in the classical case, we may conjecture that given $f\in \mathcal{T}_p^{(\alpha)}(t^\alpha)$ and $g\in \mathcal{T}_{q}^{(\alpha)}(t^\alpha)$ then $fg\in \mathcal{T}_r^{(\alpha)}(t^\alpha)$ and
    $$
\Vert fg\Vert_{\alpha,r}\le C_\alpha \Vert f\Vert_{\alpha,p}\Vert g\Vert_{\alpha,q},
$$
where $C_\alpha$ is a positive constant and $r, p, q\in [1,\infty]$ such that ${1\over r}={1\over p}+{1\over q}$. Note that for $\alpha\in \NN$, the statement holds in a straightforward way from the usual Leibniz formula and  H\"{o}lder inequality.

\item[(iii)] In Theorem \ref{convProd} we give a expression of $\mathcal{S}_{n,m} (f\otimes g)$ in terms of sums and products of different Stieltjes transform  $\mathcal{S}_{j,m} (f)$ and $\mathcal{S}_{l,m} (g)$ with $1 \leq j,l \leq m$. It is natural to conjecture that a expression of $\mathcal{S}_{\beta,\mu} (f\otimes g)$ in terms of integrals of $\mathcal{S}_{\gamma,\mu} (f)$ and $\mathcal{S}_{\nu,\mu} ( g)$ holds
with $\beta \le\gamma, \nu\leq \mu$. We will need a decomposition of the function $h_{\beta, \mu}$ in product of separate variable functions on $s$ and $u$,  where $$h_{\beta, \mu}(t,s,u)=\displaystyle{s^{\beta-1}(t+u)^\mu- u^{\beta-1}(t+s)^\mu\over {u-s}}, \qquad t>0,$$
see Lemma \ref{general}.

\end{itemize}
\subsection*{Acknowledgements}

Authors thank Aristos Siskakis and Jos\'e E. Gal\'e for several ideas, comments and usual references which have led to obtain some of these results and the final improvement of the paper.

\end{document}